\documentclass[12pt]{article}

\usepackage{amsmath,amsthm}
\usepackage{amssymb}
\usepackage{epsfig}

\usepackage[utf8x]{inputenc}
\usepackage{easybmat}
\usepackage{graphicx,hyperref}

\usepackage{url}

\newtheorem{theorem}{Theorem}[section]
\newtheorem{lemma}[theorem]{Lemma}

\theoremstyle{definition}
\newtheorem{define}[theorem]{Definition}
\newtheorem{remark}[theorem]{Remark}

\newcommand{\casefrac}[2]{ {\textstyle{ \frac{#1}{#2} } } } 

\newcommand{\pf}{\noindent {\bf Proof: }}
\newcommand{\enpf}{\hfill $\Box$ \vspace{.2in} }
\newcommand{\vs}{\vspace{0.15in}}
\newcommand{\nin}{\noindent}

\begin{document}

\title{Classifying Four-Body Convex Central Configurations}

\author{Montserrat Corbera\thanks{Departament de Tecnologies Digitals i de la Informaci\'{o}, 
Universitat de Vic, montserrat.corbera@uvic.cat}
\and
Josep M. Cors\thanks{Departament de Matem\`{a}tiques,
Universitat Polit\`{e}cnica de Catalunya, cors@epsem.upc.edu}
\and  Gareth E. Roberts\thanks{
Dept. of Mathematics and Computer Science,
College of the Holy Cross, groberts@holycross.edu}}

\maketitle

\begin{abstract}
We classify the full set of convex central configurations in the Newtonian four-body problem.
Particular attention is given to configurations possessing some type of symmetry or defining geometric property.
Special cases considered include kite, trapezoidal, co-circular, equidiagonal, orthodiagonal, and bisecting-diagonal configurations.
Good coordinates for describing the set are established.   We use them to prove that the set of four-body convex central configurations with positive masses
is three-dimensional, a graph over a domain~$D$ that is the union of elementary regions in~$\mathbb{R}^{+^3}$.
\end{abstract}

%{\bf MSC Classifications:} 70F10, 70F15, 37N05

%\vspace{.1in}


%%%%%%%%%%%%%%%%%%%%%%%%%%%%%%%%%%%%%%%%%%%%%%%%%%%%%%%%%%%%%%%%%%%%%%%%%%%%%%%%%%%%
\section{Introduction}
%%%%%%%%%%%%%%%%%%%%%%%%%%%%%%%%%%%%%%%%%%%%%%%%%%%%%%%%%%%%%%%%%%%%%%%%%%%%%%%%%%%

The study of central configurations in the Newtonian $n$-body problem is an active subfield of celestial mechanics.  
A configuration is {\em central} if the gravitational force on each body is a common scalar multiple of its position vector
with respect to the center of mass.   Perhaps the most well known example is the equilateral triangle solution of Lagrange, discovered
in 1772, consisting of three bodies of arbitrary mass located at the vertices of an equilateral triangle~\cite{lagrange}.  Released from rest,
any central configuration will collapse homothetically toward its center of mass, ending in total collision.  In fact, any solution of
the $n$-body problem containing a collision must have its colliding bodies asymptotically approaching a central configuration~\cite{saari}.
On the other hand, given the appropriate initial velocities, a central configuration can rotate rigidly about its center of mass, 
generating a periodic solution known as a {\em relative equilibrium}.   These are some of the only explicitly known solutions in the $n$-body problem.  
For more background on central configurations and their special properties, see~\cite{albouy-chen, meyer, rick-book, rick2, saari, schmidt, wintner} and the references therein.

In this paper we focus on four-body convex central configurations.  A configuration is {\em convex} if no body lies inside or on the convex hull of the other three bodies
(e.g., a rhombus or a trapezoid); otherwise, it is called {\em concave}.   Most of the results on four-body central configurations are either for a specific choice of 
masses or for a particular geometric type of configuration.  For instance, Albouy proved that all of the four-body equal mass central configurations possess a line
of symmetry.  This in turn allows for a complete solution to the equal mass case~\cite{albouy-symm1, albouy-symm2}.  
Albouy, Fu, and Sun showed that a convex central configuration with two equal masses opposite each other is symmetric,
with the equal masses equidistant from the line of symmetry~\cite{albouy}.  Recently, Fernandes, Llibre, and Mello proved that a convex central configuration
with two pairs of adjacent masses must be an isosceles trapezoid~\cite{fernand}.  
A numerical study for the number of central configurations in the four-body problem with arbitrary masses was done by Sim\'{o} in~\cite{simo}.
Other studies have focused on examples with one infinitesimal mass,
solutions of the planar restricted four-body problem~\cite{leandro1, leandro2}.

In terms of restricting the problem to a particular shape, Cors and Roberts classified the four-body co-circular central configurations in~\cite{pitu-gr} while
Corbera et al. recently studied the trapezoidal solutions~\cite{CC-Trap} (see also~\cite{santoprete}).   
Symmetric central configurations are often the easiest to analyze.  The regular $n$-gon ($n \geq 4$) is a central configuration as long as the masses are all equal.
A {\em kite} is a symmetric quadrilateral with two bodies lying on the axis of symmetry and the other two bodies positioned equidistant from it.
A kite may either be convex or concave.  In the convex case, the diagonals are always perpendicular.    
A recent investigation of the kite central configurations (both convex and concave) was carried out in~\cite{erdi}.

One of the major results in the study of convex central configurations is that they exist.  MacMillan and Bartky showed that for any four masses
and any ordering of the bodies, there exists a convex central configuration~\cite{Mac}.  This was proven again later in
a simpler way by Xia~\cite{Xia}.  It is an open question as to whether this solution must be unique.
This is problem~10 on a published list of open questions in celestial mechanics~\cite{albouy-pblms}.
Hampton showed that for any four choices of positive masses there exists a concave central configuration~\cite{hampton}.
Uniqueness does not hold in the concave setting as the example of an equilateral triangle with an arbitrary mass at the center
illustrates.  Finally, Hampton and Moeckel showed that given four positive masses, the number of equivalence
classes of central configurations under rotations, translations, and dilations is finite~\cite{rick-hamp}.

Here we study the full space of four-body convex central configurations, focusing on how various geometrically-defined
classes fit within the larger set.   We introduce simple yet effective coordinates to describe the space up to an isometry, rescaling, or
relabeling of the bodies.  Three radial coordinates $a, b$, and $c$ represent the distance from three of the bodies, respectively,  
to the intersection of the diagonals.  The remaining coordinate $\theta$ is the angle between the two diagonals.  
Positivity of the masses imposes various constraints on the coordinates.  
We find a simply connected domain $D \subset \mathbb{R}^{+^3}$, the union of four elementary regions,  
such that for any $(a,b,c) \in D$, there exists a unique angle~$\theta$ which makes the configuration central
with positive masses.  The angle $\theta = f(a,b,c)$ is a differentiable function on the interior of~$D$.  Thus the set of convex
central configurations with positive masses is the graph of a function of three variables.  We also prove that
$\pi/3 < \theta \leq \pi/2$, with $\theta = \pi/2$ if and only if the configuration is a kite.

One of the surprising features of our coordinate system is the simple linear and quadratic equations that define various
classes of quadrilaterals.  The kite configurations lie on two orthogonal planes that intersect in the
family of rhombii solutions.  These planes form a portion of the boundary of~$D$.
The co-circular and trapezoidal configurations each lie on saddles in~$D$, while the equidiagonal solutions
are located on a plane.  These three types of configurations intersect in a line corresponding to the isosceles trapezoid family.  
Our work provides a unifying structure for the set of convex central configurations and a clear picture of how the
special sub-cases are situated within the broader set.

The paper is organized as follows.  In the next section we develop the equations for a four-body central configuration using
mutual distance coordinates.   In Section~3 we introduce our coordinate system and study the important domain~$D$, proving that
$\theta$ is a differentiable function on~$D$.  We also verify the bounds on $\theta$ and show that it increases with~$c$. 
Section~4 focuses on four special cases---kite, trapezoidal, 
co-circular, and equidiagonal configurations---and how they fit together within~$D$.

Figure~\ref{Fig:Dproj} and all of the three-dimensional plots in this paper were created using Matlab~\cite{matlab}.
All other figures were made using the open-source software Sage~\cite{sage}.

%%%%%%%%%%%%%%%%%%%%%%%%%%%%%%%%%%%%%%%%%
\section{Four-Body Planar Central Configurations}
\label{sec:ccc}
%%%%%%%%%%%%%%%%%%%%%%%%%%%%%%%%%%%%%%%%%

Let $q_i \in \mathbb{R}^2$ and $m_i$ denote the position and mass, respectively, of the
$i$th body.  We will assume that $m_i > 0 \; \forall i$, while recognizing that
the zero-mass case is important for defining certain boundaries of our space.  
Let $r_{ij} = ||q_i - q_j||$ represent the distance between the $i$th and $j$th bodies.
If $M = \sum_{i=1}^n m_i$ is the sum of the masses, then the {\em center of mass} is given
by $c = \frac{1}{M} \sum_{i=1}^n m_i q_i$.  The motion of the bodies is governed by the 
Newtonian potential function 
$$
U(q) \; = \;  \sum_{i < j}^n  \; \frac{m_i m_j}{r_{ij} }.
$$
The {\em moment of inertia} with respect to the center of mass is given by
$$
I(q) \; = \;  \sum_{i=1}^{n} \; m_i \| q_i -  c \|^2 \; = \;
\frac{1}{M} \sum_{i < j}  \;  m_i m_j r_{ij}^2 .
$$
This can be interpreted as a measure of the relative size of the configuration.

There are several ways to describe a central configuration.  We follow the topological 
approach. 

\begin{define}
A planar {\em central configuration} $(q_1, \ldots, q_n) \in \mathbb{R}^{2n}$
is a critical point of $U$ subject to the constraint $I = I_0$,  where $I_0 > 0$ is a constant.  
\end{define}

It is important to note that, due to the invariance of $U$ and $I$ under isometries, any rotation, translation, or
scaling of a central configuration still results in a central configuration.

%%%%%%%%%%%%%%%%%%%
\subsection{Mutual distance coordinates}
%%%%%%%%%%%%%%%%%%%

Our derivation of the equations for a four-body central configuration follows the nice exposition of
Schmidt~\cite{schmidt}.  In the case of four bodies, the six mutual distances $r_{12}, r_{13}, r_{14}, r_{23}, r_{24}, r_{34}$ turn out to be excellent
coordinates.   They describe a configuration in the plane as long as the
Cayley-Menger determinant
$$
V \; = \;  \left|  \begin{array}{ccccc}
		0 & 1  & 1   & 1  & 1 \\[0.05in]
	         1 & 0  & r_{12}^2 & r_{13}^2 & r_{14}^2  \\[0.05in]
                  1 & r_{12}^2 &  0  & r_{23}^2  & r_{24}^2 \\[0.05in]
                  1 & r_{13}^2 & r_{23}^2 & 0 & r_{34}^2    \\[0.05in]  
                  1 & r_{14}^2 & r_{24}^2 & r_{34}^2 & 0 
              \end{array} \right|
$$
vanishes and the triangle inequality $r_{ij} + r_{jk} > r_{ik}$ holds for any choice of
indices with $i \neq j \neq k$.   The constraint $V=0$ is necessary for locating planar central configurations; 
without it, the only critical points of $U$ restricted to $I = I_0$ are regular tetrahedra (a spatial central configuration for
any choice of masses).  Therefore, we search for critical points of the function
\begin{equation}
U + \lambda(I - I_0) + \mu V
\label{eq:cps}
\end{equation}
satisfying $I = I_0$ and $V = 0$, where $\lambda$ and $\mu$ are Lagrange multipliers.

A useful formula involving the Cayley-Menger determinant is
\begin{equation}
\frac{ \partial V}{ \partial r_{ij}^2} \; = \;  - 32 \, A_i  A_j \, ,
\label{eq:cayley-deriv}
\end{equation}
where $A_i$ is the signed area of the triangle whose vertices contain all bodies
except for the $i$th body.    Formula~(\ref{eq:cayley-deriv}) holds only when restricting to planar configurations.

Differentiating~(\ref{eq:cps}) with respect to $r_{ij}$ and applying formula~(\ref{eq:cayley-deriv}) yields 
\begin{equation}
m_i m_j (s_{ij} - \lambda^{'}) \; = \;  \sigma A_i A_j,
\label{eq:cc}
\end{equation}
where $s_{ij} = r_{ij}^{-3}, \lambda^{'} = 2\lambda/M,$ and $\sigma = -64 \mu$.
Arranging the six equations of~(\ref{eq:cc}) as 
\begin{equation}
\begin{split}
& m_1 m_2(s_{12} - \lambda') = \sigma A_1 A_2, \qquad   m_3  m_4(s_{34} - \lambda') = \sigma A_3 A_4,\\
& m_1 m_3(s_{13} - \lambda') = \sigma A_1 A_3, \qquad   m_2  m_4(s_{24} - \lambda') = \sigma A_2 A_4,\\
& m_1 m_4(s_{14} - \lambda') = \sigma A_1 A_4, \qquad   m_2  m_3(s_{23} - \lambda') = \sigma A_2 A_3,
\end{split}
\label{eq:ccGroup}
\end{equation}
and multiplying together pairwise yields the well-known Dziobek relation~\cite{dzio}
\begin{equation}
(s_{12} - \lambda^{'})(s_{34} - \lambda^{'}) \; = \;  (s_{13} - \lambda^{'})(s_{24} - \lambda^{'})\;  = \;  (s_{14} - \lambda^{'})(s_{23} - \lambda^{'}).
\label{equs:dzio}
\end{equation}
This assumes that the masses and areas are all non-zero.  
Eliminating $\lambda^{'}$ from~(\ref{equs:dzio}) produces the important equation
\begin{equation}
(r_{24}^3 - r_{14}^3)(r_{13}^3 - r_{12}^3) (r_{23}^3 - r_{34}^3)  \; = \;  (r_{12}^3 - r_{14}^3) (r_{24}^3 - r_{34}^3) (r_{13}^3 - r_{23}^3) .
\label{Eq:consist}
\end{equation}

In some sense, equation~(\ref{Eq:consist}) is the defining equation for a four-body central configuration.  
It or some equivalent variation can be found in many papers and texts (e.g., see p.~278 of~\cite{wintner}.)
Equation~(\ref{Eq:consist}) is clearly necessary given the above derivation.  However, it is also sufficient assuming the six mutual
distances describe an actual configuration in the plane.  The only other restrictions required on the~$r_{ij}$ are those
that insure solutions to system~(\ref{eq:ccGroup}) yield positive masses, as explained in the next section.

%%%%%%%%%%%%%%%%%%%%%%%
\subsection{Restrictions on the mutual distances}
%%%%%%%%%%%%%%%%%%%%%%%

For the remainder of the paper we will restrict our attention to four-body convex central configurations.  We will assume the bodies are ordered
consecutively in the counterclockwise direction.  This implies that the lengths of the diagonals are $r_{13}$ and $r_{24}$, while 
the four exterior side lengths are $r_{12}, r_{23}, r_{14}$, and~$r_{34}$.   With this choice of labeling, we always have
$A_1, A_3 > 0$ and $A_2, A_4 < 0$.  We will also assume, without loss of generality,
that the largest exterior side length is $r_{12}$.

First note that $\sigma \neq 0$.  If this was not the case, then equation~(\ref{eq:cc}) and nonzero masses would imply 
that all $r_{ij}$ are equal, which is the regular tetrahedron solution.   If $\sigma < 0$, then system~(\ref{eq:ccGroup}) and positive masses implies
\begin{equation}
r_{12}, r_{14}, r_{23}, r_{34} \; < \; \frac{1}{\sqrt[3]{\lambda'}}  \; < \;  r_{13}, r_{24}  \,  .
\label{Ineq:CCdist}
\end{equation}
This means the two diagonals are strictly longer than any of the exterior sides.  On the other hand, if we assume that
$\sigma > 0$, then the inequalities in~(\ref{Ineq:CCdist}) would be reversed.  But such a configuration is impossible since
it violates geometric properties of convex quadrilaterals such as $r_{13} + r_{24} > r_{12} + r_{34}$ (see Lemma~2.3 in~\cite{HRS}).

In addition to~(\ref{Ineq:CCdist}), further restrictions on the exterior side lengths follow from the Dziobek equation
\begin{equation}
(s_{12} - \lambda^{'})(s_{34} - \lambda^{'}) \; = \;  (s_{14} - \lambda^{'})(s_{23} - \lambda^{'}).
\label{eq:sij}
\end{equation}
Since $r_{12}$ is the largest exterior side length, we  have $r_{12} \geq r_{14}$ and
$s_{14} - \lambda' \geq s_{12} - \lambda' > 0$.  It follows that $s_{34} - \lambda' \geq s_{23} - \lambda'$, otherwise
equation~(\ref{eq:sij}) is violated.  We conclude that $r_{23} \geq r_{34}$.  A similar argument shows that $r_{12} \geq r_{23}$
implies that $r_{14} \geq r_{34}$.  Hence, the shortest exterior
side is always opposite the longest one, with equality only in the case of a square.
In sum, for our particular arrangement of the four bodies, any convex central configuration with positive masses must satisfy
\begin{equation}
r_{13}, r_{24} \;  >  \;  r_{12} \; \geq \; r_{14}, r_{23} \; \geq  \;  r_{34}.
\label{geomCond}
\end{equation}

According to the Dziobek equations~(\ref{equs:dzio}), 
$$
\lambda' \; = \;  \frac{ s_{12} s_{34} - s_{13} s_{24} }{s_{12} + s_{34} - s_{13} - s_{24}}  
\; = \;  \frac{ s_{12} s_{34} - s_{14} s_{23} }{s_{12} + s_{34} - s_{14} - s_{23}}  
\; = \;  \frac{ s_{13} s_{24} - s_{14} s_{23} }{s_{13} + s_{24} - s_{14} - s_{23}} \, .
$$
These expressions generate nice formulas for the ratios between the masses.  From system~(\ref{eq:ccGroup}), a short calculation gives
\begin{equation}
\frac{m_2}{m_1} \; = \;  -\frac{A_2(s_{14} - s_{13})}{ A_1 (s_{23} - s_{24}) },  \quad
\frac{m_3}{m_1} \; = \;  \frac{A_3(s_{14} - s_{12})}{ A_1 (s_{34} - s_{23}) },  \quad 
\frac{m_4}{m_1} \; = \;  -\frac{A_4(s_{12} - s_{13})}{ A_1 (s_{34} - s_{24}) } \, 
\label{eq:massRatios1}
\end{equation}
and
\begin{equation}
\frac{m_3}{m_2} \; = \;  -\frac{A_3(s_{12} - s_{24})}{ A_2 (s_{34} - s_{13}) },  \quad
\frac{m_4}{m_2} \; = \;  \frac{A_4(s_{23} - s_{12})}{ A_2 (s_{34} - s_{14}) },  \quad 
\frac{m_4}{m_3} \; = \;  -\frac{A_4(s_{23} - s_{13})}{ A_3 (s_{14} - s_{24}) } \, .
\label{eq:massRatios2}
\end{equation}
Due to equation~(\ref{Eq:consist}), these formulas are consistent with each other.
They are all well-defined for configurations satisfying the inequalities in~(\ref{geomCond}) unless
$s_{34} = s_{23}$ (and thus $s_{12} = s_{14}$), or $s_{34} = s_{14}$ (and thus $s_{12} = s_{23}$).
For these special cases, which correspond to symmetric kite configurations, we use the alternative formulas
\begin{equation}
\frac{m_3}{m_1} \; = \;  \frac{A_3(s_{12} - s_{13})(s_{14} - s_{24})}{ A_1 (s_{23} - s_{13})(s_{34} - s_{24}) } \quad \mbox{and} \quad
\frac{m_4}{m_2} \; = \;  \frac{A_4(s_{23} - s_{13})(s_{12} - s_{24})}{ A_2 (s_{34} - s_{13})(s_{14} - s_{24}) }  \, .
\label{eq:ratioAlt}
\end{equation}

The formulas obtained for the mass ratios explain why equation~(\ref{Eq:consist}) is also sufficient for obtaining a central configuration.   
If the mutual distances $r_{ij}$ satisfy both (\ref{geomCond}) and~(\ref{Eq:consist}), then the mass ratios (which are positive), 
are given uniquely by (\ref{eq:massRatios1}),  (\ref{eq:massRatios2}), or~(\ref{eq:ratioAlt}).
We can then work backwards and check that system~(\ref{eq:ccGroup}) is satisfied so that the configuration is indeed central.

%%%%%%%%%%%%%%%%%%%%%%%%%%%%%%%%%%%%%%%%%%%%
\section{The Set of Convex Central Configurations}
\label{Sec:ConvexCC}
%%%%%%%%%%%%%%%%%%%%%%%%%%%%%%%%%%%%%%%%%%%%

We now describe the full set of convex central configurations with positive masses, showing it is three-dimensional, the graph
of a differentiable function of three variables.

%%%%%%%%%%%%%%
\subsection{Good coordinates}
%%%%%%%%%%%%%%

We begin by defining simple, but extremely useful coordinates.
Since the space of central configurations is invariant under isometries, we may apply 
a rotation and translation to place bodies 1 and~3 on the horizontal axis, with the origin located at the intersection of the two diagonals.  
It is also permissible to rescale the configuration so that $q_1 = (1,0)$.  This alters the value of the Lagrange multipliers,
but preserves the special trait of being central.

\begin{figure}[h!]
 \centering
 \includegraphics[height=10cm,keepaspectratio=true]{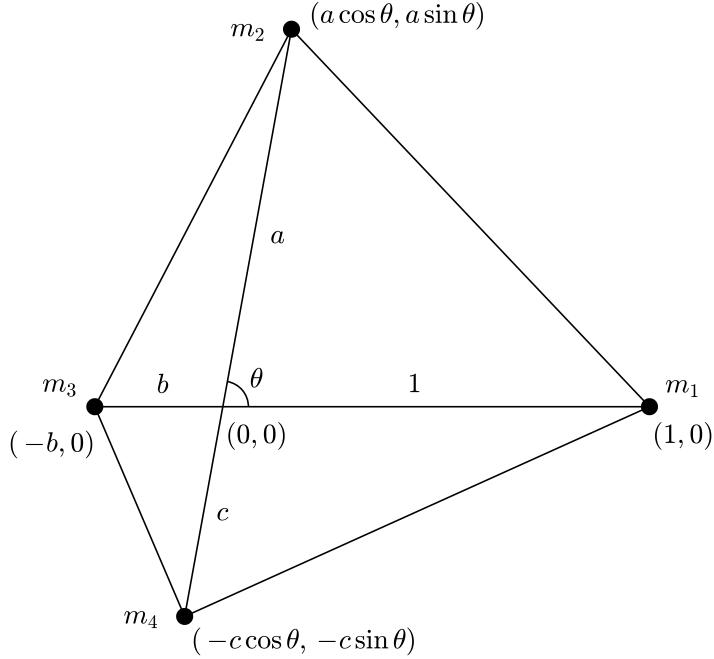}
 \vspace{-0.15in}
 \caption{Coordinates for a convex configuration of four bodies:  three radial variables $a, b, c > 0$ and an angular variable $\theta \in (0,\pi)$.}
 \label{Fig:setup} 
\end{figure}

Define the remaining three bodies to have positions $q_2 = (a \cos \theta, a \sin \theta), q_3 = (-b, 0)$, and $q_4 = (-c \cos \theta, -c \sin \theta)$, 
where $a, b, c$ are radial variables and $\theta \in (0, \pi)$ is an angular variable (see Figure~\ref{Fig:setup}).
If one or more of the three radial variables is negative, then the configuration becomes concave or the ordering of the bodies changes.
If one or more of the radial variables vanish, then the configuration contains a subset that is collinear or some bodies coalesce (e.g.,
$b=c=0$ implies $r_{34} = 0$).    Thus, we will assume throughout the paper that $a > 0, b > 0,$ and $c > 0$.
The coordinates $(a,b,c,\theta)$ turn out to be remarkably well-suited for describing different classes of quadrilaterals 
that are also central configurations (see Section~\ref{Sec:Shapes}).

\vspace{-0.15in}

%----------------------------------------------------
%------------------------------------------------

In our coordinates, the six mutual distances are given by
\begin{eqnarray}
r_{12}^2 = a^2 -2a \cos \theta + 1,    \quad &  r_{23}^2 = a^2 + 2ab \cos \theta + b^2, &  \quad  r_{13} = b + 1,  \label{eq:MutDist1} \\[0.07in]
r_{14}^2 = c^2 + 2c \cos \theta + 1,   \quad &  r_{34}^2 = b^2 - 2bc \cos \theta + c^2,  &  \quad  r_{24} = a + c.  \label{eq:MutDist2}
\end{eqnarray}
Based on equation~(\ref{Eq:consist}), define $F$ to be the function
$$
F(a,b,c,\theta) \; = \;  (r_{24}^3 - r_{14}^3)(r_{13}^3 - r_{12}^3)(r_{23}^3 - r_{34}^3) -  (r_{12}^3 - r_{14}^3)(r_{24}^3 - r_{34}^3)(r_{13}^3 - r_{23}^3),
$$
where each mutual distance is treated as a function of the variables $a, b, c,$ and $\theta$.

The previous discussion justifies the following lemma.
\begin{lemma}
Let ${\cal C}$ and $E$ denote the sets
\begin{eqnarray*}
{\cal C} & = &  \{  (a,b,c,\theta) \in \mathbb{R}^{+^3} \! \times (0, \pi): r_{13}, r_{24}   >   r_{12}  \geq r_{14}, r_{23}  \geq r_{34}   \} ,  \\
E  & = &  \{ s = (a,b,c,\theta) \in \mathbb{R}^{+^3} \! \times (0, \pi): s \in {\cal C} \mbox{ and } F(s) = 0 \}.
\end{eqnarray*}
Any point in $E$ corresponds to a four-body convex central configuration with positive masses.  Moreover, up to an isometry, rescaling, or relabeling of the bodies,
$E$ contains all such configurations.
\end{lemma}

%%%%%%%%%%%%%%
\subsection{Defining the domain $D$}
%%%%%%%%%%%%%%

We will find a set $D \subset \mathbb{R}^{+^3}$ such that for each $(a,b,c) \in D$, there exists a unique angle~$\theta$ which
makes the configuration central.  Specifically, we prove that there exists a differentiable function $\theta = f(a, b, c)$ with domain $D$, whose
graph is equivalent to $E$.  In order to define~$D$, we use the mutual distance inequalities in~(\ref{geomCond}) to eliminate the angular variable~$\theta$.

\begin{lemma}
The inequalities in~(\ref{geomCond}) imply the following conditions on the positive variables $a, b, c$:
\begin{eqnarray}
r_{12} \; \geq \;  r_{14} \; \mbox{ and } \; r_{23} \; \geq \;  r_{34} & \quad \Longrightarrow \quad &  a \; \geq \; c ,  \label{Cond1} \\[0.1in]
r_{12} \; \geq \;  r_{23} \; \mbox{ and } \; r_{14} \; \geq \;  r_{34} & \quad \Longrightarrow \quad &  b \; \leq \; 1 , \label{Cond2} \\[0.1in]
r_{13}  \; >  \; r_{12} \;  \geq   \;  r_{14}   & \quad \Longrightarrow \quad &  c \; < \;  \frac{1}{a} (b^2 + 2b), \label{Cond3}\\[0.1in]
r_{13}  \; >  \; r_{12} \;  \geq   \;  r_{23}  \; \mbox{ and } \; a > 1 & \quad \Longrightarrow \quad &  b \; >  \; \frac{1}{2}(-1 + \sqrt{4a^2 - 3}\, ),  \label{Cond5} \\[0.1in] 
r_{24} \;  >  \; r_{12} \;  \geq   \;  r_{14}  \; \mbox{ and } \;  0 < a < 1 & \quad \Longrightarrow \quad &   c \; > \; \frac{1}{2} (-a + \sqrt{4 - 3a^2}\,),  \label{Cond4} \\[0.1in]
r_{24} \; >  \;  r_{12} \;  \geq   \;  r_{23}   & \quad \Longrightarrow \quad &   c \; > \; -a + \sqrt{a^2 + b} \, .  \label{Cond6}
\end{eqnarray}
\label{Lemma:DistBds}
\end{lemma}

\pf
From equations (\ref{eq:MutDist1}) and~(\ref{eq:MutDist2}) we compute that
\begin{eqnarray}
r_{12}^2 - r_{14}^2 \; = \; (a + c)(a - c - 2 \cos \theta),    &  &  r_{12}^2 - r_{23}^2 \; = \; (1 + b)(1 - b - 2a \cos \theta), \label{ident1} \\[0.1in]
r_{23}^2 - r_{34}^2 \; = \;  (a + c)(a - c + 2b \cos \theta), &  &  r_{14}^2 - r_{34}^2 \; = \; (1 + b)(1 - b + 2c \cos \theta).  \label{ident2}
\end{eqnarray}
Since $a, b,$ and $c$ are all positive, $r_{12} \geq r_{14}$ and $r_{23} \geq r_{34}$ together imply that
\begin{equation}
a - c \;  \geq \; \mbox{max}\{2 \cos \theta, -2b \cos \theta\} \; \geq \; 0 .
\label{MaxBd1}
\end{equation}
Similarly, $r_{12} \geq r_{23}$ and $r_{14} \geq r_{34}$ imply
\begin{equation}
1 - b \; \geq \;  \mbox{max}\{2 a \cos \theta, -2c \cos \theta\} \; \geq \; 0.
\label{MaxBd2}
\end{equation}
This proves implications (\ref{Cond1}) and~(\ref{Cond2}).

Next, equations (\ref{eq:MutDist1}) and~(\ref{eq:MutDist2}) yield
\begin{equation}
r_{13}^2 - r_{12}^2 \; = \;  b^2 + 2b - a^2 + 2a \cos \theta \quad \mbox{and} \quad
r_{24}^2 - r_{12}^2 \; = \;  c^2 + 2ac - 1 + 2a \cos \theta .
\label{ident3}
\end{equation}
Since $r_{12} \geq r_{14}$, the first equation in~(\ref{ident1}) gives $a - 2 \cos \theta \geq c$ or $a^2 - 2a \cos \theta \geq ac$. Then
$r_{13} > r_{12}$ implies that 
\begin{equation}
b^2 + 2b \; > \;  a^2 - 2a \cos \theta \; \geq \; ac,
\label{ineq:r13r12r14}
\end{equation}
which verifies~(\ref{Cond3}).

Similarly, $r_{12} \geq r_{23}$ and the second equation in~(\ref{ident1}) yields $-2a \cos \theta \geq b - 1$.  Then 
$r_{13} > r_{12}$ implies that 
\begin{equation}
b^2 + 2b - a^2 \;  > \;  -2a \cos \theta \; \geq \;  b - 1,
\label{ineq:r13r12r23}
\end{equation}
which yields
\begin{equation}
b^2 + b + 1 - a^2 \; > \; 0.
\label{ineq:Quad2}
\end{equation}
Since $b$ and $a$ are both positive, inequality~(\ref{ineq:Quad2}) clearly holds if $a \leq 1$.  However, for any fixed choice of $a > 1$, the value
of $b$ must be chosen strictly greater than the largest root of the quadratic $Q_a(b) = b^2 + b + 1 - a^2$.  This root is
$\frac{1}{2}(-1 + \sqrt{4a^2 - 3})$, which verifies implication~(\ref{Cond5}).

Next, $r_{24} > r_{12} \geq r_{14}$ yields
\begin{equation}
c^2 + 2ac -1 + a^2 \; > \;  - 2a \cos \theta + a^2 \;  \geq \;  ac,
\label{ineq:r24r12r14}
\end{equation}
which in turn gives
\begin{equation}
c^2 + ac + a^2 - 1 \; > \; 0.
\label{ineq:Quad}
\end{equation}
Since $a$ and $c$ are both positive, inequality~(\ref{ineq:Quad}) clearly holds if $a \geq 1$.  However, for any fixed choice of $a \in (0,1)$, the value
of $c$ must be chosen strictly greater than the largest root of the quadratic $Q_a(c) = c^2 + ac + a^2 - 1$.  This root is
$\frac{1}{2}(-a + \sqrt{4 - 3a^2})$, which proves~(\ref{Cond4}).

Finally, $r_{24} >  r_{12} \geq r_{23}$ implies that 
\begin{equation}
c^2 + 2ac -1 \; > \;  -2a \cos \theta \; \geq \;  b - 1,
\label{ineq:r24r12r23}
\end{equation}
which gives
\begin{equation}
c^2 + 2ac - b \; > \; 0.
\label{ineq:Quad3}
\end{equation}
Since $b > 0$, $c$ must be chosen greater than the largest root of the quadratic $Q_{a,b}(c) = c^2 + 2ac - b$.  This root is
$-a + \sqrt{a^2 + b}$, which verifies~(\ref{Cond6}) and completes the proof.
\enpf

The combined inequalities between the radial variables $a, b,$ and $c$ given in (\ref{Cond1}) through~(\ref{Cond6}), along with
$a > 0, b > 0,$ and $c > 0$, define a bounded set~$D \subset \mathbb{R}^{+^3}$.  We will show that this set is the
domain of the function $\theta = f(a,b,c)$ and the projection of $E$ into $abc$-space.

\vs

\fbox{\; \;  \parbox{6.4in}{ 
\begin{define}
Let $D = D_1 \cup D_2$ denote the three-dimensional region, where
\begin{eqnarray*}
D_1 & = & \Bigl\{ (a,b,c) \in \mathbb{R}^{+^3}: 0 < c \leq a, \, a \leq 1, \, 0 < b \leq 1,  \\[0.1in]
&  &    \frac{1}{2} (-a + \sqrt{4 - 3a^2}\,) < c <  \frac{1}{a} (b^2 + 2b), \,  c > -a + \sqrt{a^2 + b} \, \Bigr\} ,  \\[0.15in]
D_2 &  = &   \Bigl\{ (a,b,c) \in \mathbb{R}^{+^3}: 0 < c \leq a, \, a > 1, \,  0 < b \leq 1, \,  c <  \frac{1}{a} (b^2 + 2b), \\[0.1in]
&   &   b > \frac{1}{2}(-1 + \sqrt{4a^2 - 3}\, ), \, c > -a + \sqrt{a^2 + b} \, \Bigr\}.
\end{eqnarray*}
\label{def:D}
\end{define}
\vspace{-0.25in}
}}

\vs

Note that $D$ is simply connected.  
Using inequalities (\ref{ineq:Quad2}), (\ref{ineq:Quad}), $c \leq a$, and $b \leq 1$, it is straight-forward to check that $D$ is contained within the box
$$
\frac{1}{\sqrt{3}} \leq a \leq \sqrt{3}, \;  0 \leq b \leq 1,  \; 0 \leq c \leq \sqrt{3} \,.
$$

\begin{figure}[t]
 \centering
 \includegraphics[height=11cm,keepaspectratio=true]{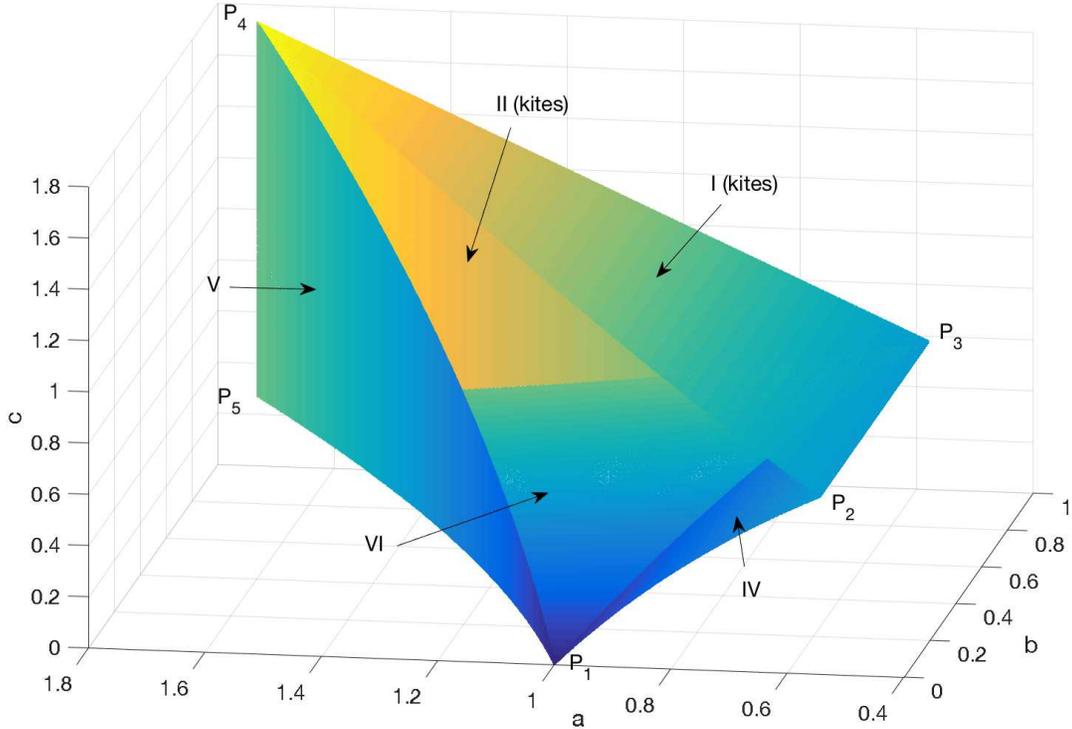}
 \vspace{-0.25in}
 \caption{The faces and vertices of $\overline{D}$ (face~III not shown to improve the perspective).   Faces II and~V are vertical.  For each point $(a, b, c) \in D$, there exists a unique angle $\theta$ that makes the corresponding
 configuration central.}
 \label{Fig:Dplot} 
\end{figure}

Let $\overline{D}$ denote the closure of~$D$.
A plot of the boundary of $D$ is shown in Figure~\ref{Fig:Dplot}.  It contains five vertices, six faces, and nine edges (six curved, three straight), 
in accordance with Poincar\'{e}'s generalization of Euler's formula $\overline{V} -\overline{E} + \overline{F} = 2$.  
The vertices of $\overline{D}$ are 
$$
\begin{array}{l}
P_1 = (1,0,0), \quad  P_2 = (\frac{1}{\sqrt{3}}, \frac{2 - \sqrt{3}}{\sqrt{3}}, \frac{1}{\sqrt{3}}), \quad P_3 = (\frac{1}{\sqrt{3}}, 1, \frac{1}{\sqrt{3}}), \\[0.1in]
P_4 = (\sqrt{3}, 1 , \sqrt{3})  , \mbox{ and } P_5 =  (\sqrt{3}, 1, 2 - \sqrt{3}),
\end{array}
$$
each of which corresponds to a symmetric central configuration with at least two zero masses.  $P_3$ and~$P_4$ are rhombii with one diagonal congruent to the common side length, while
$P_2$ and~$P_5$ are kites with horizontal and vertical axes of symmetry, respectively.  The point $P_1$ corresponds to an equilateral triangle with bodies 3 and~4 sharing a common vertex.

%----------------------------------------------------
%------------------------------------------------

%%%%%%%%%%%%%%%%%%%%%%%
\subsection{Configurations on the boundary of $D$}
\label{SubSec:BoundaryConfigs}
%%%%%%%%%%%%%%%%%%%%%%%

We now focus on points lying on the boundary of~$D$.  The next lemma shows that these points correspond to
configurations where two or more of the mutual distance inequalities in~(\ref{geomCond}) become equalities.  Moreover,
the only points for which this is true lie on the boundary of~$D$.

\begin{lemma}
Suppose that $(a,b,c,\theta)$ are chosen so that $r_{13}, r_{24}   \geq   r_{12}  \geq r_{14}, r_{23}  \geq r_{34}$ with $a \geq 1/\sqrt{3}$, $b \geq 0$, and
$c \geq 0$.  Then
\begin{eqnarray}
r_{12} = r_{14} \; \mbox{ and } \; r_{23} = r_{34} & \quad \mbox{if and only if} \quad &  a \; = \; c ,  \label{EqualCond1} \\[0.1in]
r_{12} = r_{23} \; \mbox{ and } \;  r_{14} =  r_{34} & \quad  \mbox{if and only if} \quad &  b \; = \; 1 , \label{EqualCond2} \\[0.1in]
r_{13} \; = \;  r_{12} \;  = \;  r_{14}   & \quad   \mbox{if and only if} \quad &  c \; = \;  \frac{1}{a} (b^2 + 2b), \label{EqualCond3}\\[0.1in]
r_{24} \; = \;  r_{12} \;  = \;  r_{14}   & \quad   \mbox{if and only if} \quad &   c \; = \; \frac{1}{2} (-a + \sqrt{4 - 3a^2}\,),  \label{EqualCond4} \\[0.1in]
r_{13} \; = \;  r_{12} \;  = \;  r_{23}   & \quad   \mbox{if and only if} \quad &  b \; =  \; \frac{1}{2}(-1 + \sqrt{4a^2 - 3}\, ),  \label{EqualCond5} \\[0.1in] 
r_{24} \; = \;  r_{12} \;  = \;  r_{23}   & \quad   \mbox{if and only if} \quad &   c \; = \; -a + \sqrt{a^2 + b} \, .  \label{EqualCond6}
\end{eqnarray}
\label{Lemma:Boundary}
\end{lemma}

\pf
We first note that under the assumptions of the lemma, the inequalities on $a, b,$ and $c$ from Lemma~\ref{Lemma:DistBds} are still valid, except that
the inequalities are no longer strict.

If $r_{12} = r_{14}$ and $r_{23} = r_{34}$, then the left-hand equations in (\ref{ident1}) and~(\ref{ident2}) imply $a - c = 2 \cos \theta$ and $a - c = -2b \cos \theta$, respectively.
This yields $(1+b) \cos \theta = 0$ from which it follows that  $\cos \theta = 0$ and $a = c$.  Conversely, if $a = c$, (\ref{MaxBd1}) implies that either $\cos \theta = 0$ or $b = 0$.
In the former case, $\theta = \pi/2$ and then $r_{12} = r_{14}$ and $r_{23} = r_{34}$ follows quickly.
In the latter case, inequality~(\ref{Cond3}) and $a = c$ implies that $a = c = 0$, which contradicts $a \geq 1/\sqrt{3}$.  Thus $b > 0$ and 
$r_{12} = r_{14}$ and $r_{23} = r_{34}$, proving~(\ref{EqualCond1}).

If $r_{12} = r_{23}$ and $r_{14} = r_{34}$, then the right-hand equations in (\ref{ident1}) and~(\ref{ident2}) imply $1 - b = 2 a\cos \theta$ and $1 - b = -2c \cos \theta$, respectively.
Thus, $(a+c) \cos \theta = 0$.  Since $a \geq 1/\sqrt{3}$ and $c \geq 0$, we must have $\cos \theta = 0$ and hence $b = 1$.
Conversely, if $b = 1$, (\ref{MaxBd2}) implies that either $\cos \theta = 0$, or $\cos \theta < 0$ and $c=0$.
In the former case, $\theta = \pi/2$ and then $r_{12} = r_{23}$ and $r_{14} = r_{34}$ follows quickly.
The latter case is impossible, since $c = 0$ and $b = 1$ contradicts inequality~(\ref{Cond6}).
This proves~(\ref{EqualCond2}).

If $r_{13} = r_{12}$, then the left-hand equation in~(\ref{ident3}) gives $a - 2 \cos \theta = \frac{1}{a}(b^2 + 2b)$.   Likewise, if $r_{12} = r_{14}$, then $a - 2 \cos \theta = c$.
Thus $r_{13} = r_{12} = r_{14}$ implies $c = \frac{1}{a}(b^2 + 2b)$.  Conversely, if $ac = b^2 + 2b$, then both inequalities in~(\ref{ineq:r13r12r14}) become equalities.
From this we deduce that $r_{13} = r_{12} = r_{14}$, which verifies~(\ref{EqualCond3}).

If $r_{24} = r_{12}$, then the right-hand equation in~(\ref{ident3}) gives $a(c + 2 \cos \theta) = 1 - c^2 - ac$.   Likewise, if $r_{12} = r_{14}$, then $c + 2 \cos \theta = a$.
Thus $r_{24} = r_{12} = r_{14}$ implies $c^2 + ac + a^2 - 1 = 0$.   The quadratic $Q_a(c) = c^2 + ac + a^2 - 1$ has real roots for $1/\sqrt{3} \leq a \leq 2/\sqrt{3}$, but
the smaller root is always negative for these $a$-values.  Thus $c$ must be taken to be the larger root of $Q_a(c)$.  
Conversely, if $c = \frac{1}{2} (-a + \sqrt{4 - 3a^2}\,)$, then $c^2 + ac + a^2 - 1 = 0$ and both inequalities in~(\ref{ineq:r24r12r14}) become equalities.
From this we deduce that $r_{24} = r_{12} = r_{14}$, which verifies~(\ref{EqualCond4}).  
The proof of (\ref{EqualCond5}) and~(\ref{EqualCond6}) follows in a similar fashion, using inequalities (\ref{ineq:r13r12r23}) and~(\ref{ineq:r24r12r23}), respectively.
\enpf

Lemma~\ref{Lemma:Boundary} shows that the six faces on the boundary of~$D$, labeled I through VI,  are given by the
six equations (\ref{EqualCond1}) through~(\ref{EqualCond6}), respectively.  The first two faces are the only ones belonging to~$D$ (positive masses)
and contain all of the kite configurations, where $\theta = \pi/2$.   
The remaining four faces on the boundary of~$D$
correspond to cases with one or three zero masses (see Table~\ref{Table:Faces}).   Points on these faces are interpreted as
limiting solutions of a sequence of central configurations with positive masses.  The mass values shown in Table~\ref{Table:Faces}
follow from formulas (\ref{eq:massRatios1}), (\ref{eq:massRatios2}), and~(\ref{eq:ratioAlt}). 
Here we assume that the limiting solution lies in the interior of the given face.

\renewcommand{\arraystretch}{2.00}
\begin{table}[t]
\begin{center}
\begin{tabular}{||c|c|c|c|c||}
\hline
\makebox[0.6in] {{\bf Face}} &  \makebox[1.0in]{{\bf Equation}} & \makebox[1in]{{\bf Mutual Distances}} &\makebox[1.1in]{{\bf Masses}} &\makebox[0.85in]{{\bf Vertices}} \\
\hline\hline
I &  $c = a$ & $r_{12} = r_{14}$ and $r_{23} = r_{34}$  &   $m_2 = m_4$ &  $P_2, P_3, P_4$ \\
\hline
II  &  $b = 1$ & $r_{12} = r_{23}$ and $r_{14} = r_{34}$  &  $m_1 = m_3$ &  $P_3, P_4, P_5$ \\
\hline
III  &  $c =   \frac{1}{a} (b^2 + 2b)$ & $r_{13} = r_{12} = r_{14}$  & $m_2 = m_3 = m_4 = 0$  & $P_1, P_2, P_4$  \\
\hline
IV &  $c = \frac{1}{2} (-a + \sqrt{4 - 3a^2}\,)$  & $r_{24} = r_{12} = r_{14}$  &  $m_3 = 0$  & $P_1, P_2, P_3$  \\
\hline
V  &   $b =  \frac{1}{2}(-1 + \sqrt{4a^2 - 3}\, )$ & $r_{13} = r_{12} = r_{23}$   &  $m_4 = 0$  &  $P_1, P_4, P_5$  \\
\hline
VI  &   $c = -a + \sqrt{a^2 + b}$   & $r_{24} = r_{12} = r_{23}$  &  $m_1 = m_3 = m_4 = 0$  &  $P_1, P_3, P_5$ \\
\hline\hline
\end{tabular}
\caption{The six faces on the boundary of~$D$ along with their key attributes.  Each point on the boundary has a unique angle~$\theta$ that makes
the configuration central.  On faces I and~II, $\theta = \pi/2$ (kites).  On faces III and~IV, $\theta = \cos^{-1}(\frac{a-c}{2})$, while on
faces V and~VI, $\theta = \cos^{-1}(\frac{1-b}{2a})$.}
\label{Table:Faces}
\end{center}
\end{table}
\renewcommand{\arraystretch}{1.0}

For example, suppose there is a sequence of points $x^\epsilon = (a^\epsilon, b^\epsilon, c^\epsilon)$
in the interior of $D$ converging to a point $\overline{x} = (\overline{a}, \overline{b}, \overline{c})$ located on the interior of face~V.
%($1 < \overline{a} < \sqrt{3}, \overline{b} = \frac{1}{2}(-1 + \sqrt{4\overline{a}^2-3})$, and
%$-\overline{a} + \sqrt{\overline{a}^2 + \overline{b}} < \overline{c} < \frac{1}{\overline{a}} (\overline{b}^2 + 2\overline{b})$).
This corresponds to a sequence of central configurations, each with positive masses, that limits on a configuration with
$r_{13} = r_{12} = r_{23}$.  Since $\overline{x}$ does not lie on any of the other faces on the boundary of~$D$,
the other three limiting mutual distances, $r_{24}, r_{14},$ and $r_{34}$, must be distinct from $r_{13}$ and each other.  
Moreover, the limiting values of the areas $A_i$ do not vanish because $\overline{a}, \overline{b},$ and $\overline{c}$ are
all strictly positive.  Using either (\ref{eq:massRatios1}) or~(\ref{eq:massRatios2}), it follows that the limiting mass value for $m_4$
must vanish, while the other limiting mass values are strictly positive.  A similar argument applied to the other faces determines which masses must
vanish in the limit.

Configurations on face IV or V, respectively, correspond
to equilibria of the planar, circular, restricted four-body problem with infinitesimal mass $m_3$ or $m_4$, respectively~\cite{leandro1, leandro2, kule}.
Configurations on face III or VI, respectively, correspond to relative equilibria of the $(1+3)$-body problem, where a central mass
(body 1 or 2, respectively) is equidistant from three infinitesimal masses~\cite{CC-Bifur, hall, rick1}.   Note that we have not made any assumptions
on the relative size of the masses.
Each of the six faces satisfies either $r_{12} = r_{14}$ or $r_{12} = r_{23}$.  Using identity~(\ref{ident1}), it follows that 
there is a unique value of~$\theta$ for each point on the boundary of~$D$, $\theta = \cos^{-1}(\frac{a-c}{2})$ if $r_{12} = r_{14}$
or $\theta = \cos^{-1}(\frac{1-b}{2a})$ if $r_{12} = r_{23}$.

The masses at the vertices of $\overline{D}$ are not well-defined because there are more options for the path of a limiting sequence.
For example, the point $P_4$ represents a rhombus with one diagonal ($r_{13}$) congruent to all of the exterior sides.  
Approaching $P_4$ along the line $(a, 1, a)$ as $a \rightarrow \sqrt{3}$ (a sequence of rhombii 
central configurations) yields the limiting mass values $m_2 = m_4 = 0$ and $m_1 = m_3 \neq 0$.  On the other hand, it is possible to construct
a sequence of kite central configurations on face~I with masses $m_1 = 1, m_2 = m_4 = \epsilon^2,$ and $m_3 = \epsilon$ that limits
on $P_4$ as $\epsilon \rightarrow 0$.  The first sequence has two limiting mass values that vanish while the second sequence has three.
The difference occurs because the mass ratio $m_3/m_1$ at $P_4$ is undefined in either formula (\ref{eq:massRatios1}) or~(\ref{eq:ratioAlt}).

Regardless of the particular limiting sequence, all five vertices of $\overline{D}$ will have at least two mass values that vanish in the limit.
For $P_1$, this follows from Proposition~2 in~\cite{rick-cluster}.  For the other four vertices, this fact is a consequence of
formulas (\ref{eq:massRatios1}) and~(\ref{eq:massRatios2}).  In general, note that a limiting sequence with precisely two zero masses
can only occur at vertices $P_1, P_3,$ or $P_4$.   This somewhat surprising restriction is a consequence of 
Propositions 3 and~4 in~\cite{rick-cluster} 
and the fact that the non-collinear critical points of the restricted three-body problem must form an equilateral triangle with the non-trivial masses.

\subsection{The projection of $\overline{D}$ onto the $ab$-plane }
%%%%%%%%%%%%%%%%%%%%%

\begin{figure}[htb]
 \centering
 \includegraphics[height=9.9cm,keepaspectratio=true]{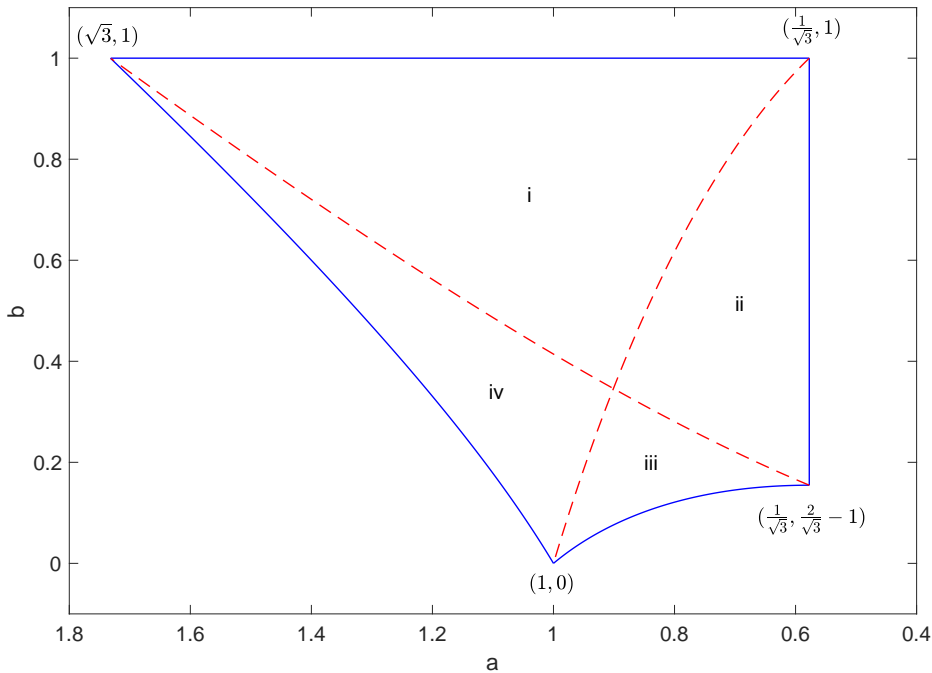}
 \vspace{-0.1in}
 \caption{The projection of $\overline{D}$ into the $ab$-plane.  
 The dashed red curves divide the region into four sub-regions over which
 $c$ is bounded by functions of $a$ and $b$.  The orientation of the $a$-axis has been reversed to match Figure~\ref{Fig:Dplot}.}
 \label{Fig:Dproj} 
\end{figure}

The set~$\overline{D}$ can be written as the union of four elementary regions in $abc$-space, where $c$ is bounded by functions of $a$ and $b$. The projection of $\overline{D}$ onto the $ab$-plane is shown
in Figure~\ref{Fig:Dproj}.  It is determined by $\frac{1}{\sqrt{3}} \leq a \leq \sqrt{3} \,$ and $l(a) \leq b \leq 1$, where $l(a)$ is the piecewise function 
$$
l(a) \; = \;  \left\{ \begin{array}{cc}
               l_1(a)  & \mbox{if $\casefrac{1}{\sqrt{3}} \leq a \leq 1$} \\[0.07in]
               l_2(a)   & \mbox{if $1 \leq a \leq \sqrt{3} \, $}.
              \end{array}
        \right.
$$
Here, $l_1(a) =  -1 + \frac{1}{2} (a + \sqrt{4 - 3a^2} \, ) $ is the projection of the intersection between faces III and~IV, and
$l_2(a) =  \frac{1}{2} (-1 + \sqrt{4a^2 - 3}\, )$ is the projection of the vertical face~V.   The edge $a = \frac{1}{\sqrt{3}}$ is the projection
of the intersection between faces I and~IV, while the edge $b=1$ is the projection of the vertical face~II.  

\vspace{-0.15in}

%----------------------------------------------------
%------------------------------------------------

The decreasing dashed curve in Figure~\ref{Fig:Dproj} is the projection of the
intersection of faces I and~III, given by $b = -1 + \sqrt{1 + a^2} \, , \frac{1}{\sqrt{3}} \leq a \leq \sqrt{3}$.
The increasing dashed curve is the projection of the
intersection of faces IV and~VI, given by
$b = 1 - \frac{3}{2}a^2 + \frac{a}{2} \sqrt{4 - 3a^2} \, ,  \frac{1}{\sqrt{3}} \leq a \leq 1$. 
These curves divide \\[0.05in]
the projection into four sub-regions over which
$c$ is bounded by different functions of $a$ and $b$, as indicated below:
$$
\begin{array}{cl}
{\rm i}   &   \quad  -a + \sqrt{a^2 + b} \; \leq \; c  \; \leq \;  a  \\[0.1in]
{\rm ii}  &   \quad  \casefrac{1}{2} (-a + \sqrt{4 - 3a^2}\,) \;  \leq \;  c  \;  \leq \;  a  \\[0.1in]
{\rm iii}  &  \quad \casefrac{1}{2} (-a + \sqrt{4 - 3a^2}\,)  \;  \leq \;  c  \;  \leq   \;  \casefrac{1}{a} (b^2 + 2b) \\[0.1in]
{\rm iv}  &  \quad   -a + \sqrt{a^2 + b}  \;  \leq  \;  c  \; \leq   \;  \casefrac{1}{a} (b^2 + 2b) \, .
\end{array}
$$

%%%%%%%%%%%%%%%%%%%%%
\subsection{$E$ is a graph $\theta = f(a,b,c)$ over $D$}
%%%%%%%%%%%%%%%%%%%%%

We now prove our main result, showing that for each $(a,b,c) \in D$, there exists a unique angle $\theta$ that makes
the configuration central.  In general, for any point $(a, b, c)$ in the interior of $D$, there is an interval of possible angles $\theta$ for which 
the mutual distance inequalities~(\ref{geomCond}) hold.  According to the identities given in (\ref{ident1}), (\ref{ident2}), and~(\ref{ident3}),
$\theta$ must be chosen to satisfy
\begin{equation}
\mbox{max} \Bigl\{ \frac{c-a}{2b}, \;  \frac{b - 1}{2c}, \;  \frac{a^2 - b^2 - 2b}{2a}, \; \frac{1 - c^2 - 2ac}{2a} \Bigl\}
\; \leq \;  \cos \theta \;  \leq \;
\mbox{min} \Bigl\{ \frac{a-c}{2}, \; \frac{1 - b}{2a} \Bigl\}
\label{ThetaBounds}
\end{equation}
in order for~(\ref{geomCond}) to be true.
The following lemma shows that condition~(\ref{ThetaBounds}) is not vacuous on the interior of $D$.

\begin{lemma}
For any point $(a,b,c)$ in the interior of $D$, define the constants $k_1$ and $k_2$ by
\begin{eqnarray*}
k_1  & = &  \max \Bigl\{ \frac{c-a}{2b}, \;  \frac{b - 1}{2c}, \;  \frac{a^2 - b^2 - 2b}{2a}, \; \frac{1 - c^2 - 2ac}{2a} \Bigl\}, \\[0.07in] 
k_2  & = &  \min \Bigl\{ \frac{a-c}{2}, \; \frac{1 - b}{2a} \Bigl\}.
\end{eqnarray*}
Then $-1 < k_1 < k_2 < 1$.
\label{lemma:AngleBounds}
\end{lemma}

\pf
On the interior of $D$ the first two quantities in the definition of $k_1$ are strictly negative while the two quantities defining $k_2$ are
strictly positive.  The inequality $(a^2 - b^2 - 2b)/(2a) < (a-c)/2$ follows from $c < (b^2 + 2b)/a$.
The inequality $(a^2 - b^2 - 2b)/(2a) < (1-b)/(2a)$ is equivalent to $b^2 + b + 1 - a^2 > 0$, which is clearly
valid for $a \leq 1$.  It also holds for $a > 1$ because $b > \frac{1}{2}(-1 + \sqrt{4a^2 - 3} \, )$.
Likewise, $(1 - c^2 - 2ac)/(2a) < (a-c)/2$ is equivalent to $c^2 + ac + a^2 - 1 > 0$, which is clearly
satisfied for $a \geq 1$.  It also holds for $0 < a < 1$ since $c > \frac{1}{2}(-a + \sqrt{4 - 3a^2} \,)$.
Finally, $(1 - c^2 - 2ac)/(2a) < (1-b)/(2a)$ is satisfied because $c > -a + \sqrt{a^2 + b}$.  This verifies that $k_1 < k_2$.

Since $a < \sqrt{3} <  2 + c$ and $1 < \frac{2}{\sqrt{3}}  < 2a + b$ on the interior of~$D$, we see that $k_2 < 1$.
Finally, $(1 - c^2 - 2ac)/(2a) > -1$ holds if $c < 1$.  But if $c \geq 1$, then $b > 0 > 1 - 2c$ implies that
$(b-1)/(2c) > -1$.  Thus, at least one of the quantities in the definition for $k_1$ is larger than~$-1$.  This shows
that $k_1 > -1$.
\enpf

Lemma~\ref{lemma:AngleBounds} shows that for any point $(a, b, c)$ in the interior of $D$, there
is an interval of $\theta$-values for which~(\ref{geomCond}) holds.   More specifically, if we let
$\theta_l = \cos^{-1} (k_2)$ and $\theta_u = \cos^{-1} (k_1)$, with $k_1, k_2$ defined as in Lemma~\ref{lemma:AngleBounds},
then for any $\theta \in (\theta_l, \theta_u)$, we have $(a,b,c,\theta) \in {\cal C}$.

Recall that
$$
F(a,b,c,\theta) \; = \;  (r_{24}^3 - r_{14}^3)(r_{13}^3 - r_{12}^3)(r_{23}^3 - r_{34}^3) -  (r_{12}^3 - r_{14}^3)(r_{24}^3 - r_{34}^3)(r_{13}^3 - r_{23}^3),
$$
and that $E = \{ s = (a,b,c,\theta) \in \mathbb{R}^{+^3} \! \times (0, \pi): s \in {\cal C} \mbox{ and } F(s) = 0 \}$ represents the set of 
convex central configurations with positive masses.

\begin{theorem}
Suppose that $(a, b, c) \in D$.  Then there exists a unique angle $\theta$ such that $(a, b, c, \theta)$ determines
a central configuration.  More precisely, the set of four-body convex central configurations with positive masses is the graph of a differentiable
function $\theta = f(a,b,c)$.  The domain of this function is~$D$, which is equivalent to the projection of $E$ onto $abc$-space. 
\label{Thm:Main}
\end{theorem}

\pf
Fix a point $(a,b,c)$ in the interior of $D$ and treat $F = F(\theta)$ as a one-variable function.  
We will show that $F$ has a unique root~$\theta$ satisfying the inequalities in~(\ref{ThetaBounds}).

\vs

\nin  {\bf (i) Existence:}
Suppose that $\theta$ is taken to be $\theta_l = \cos^{-1} (k_2)$.
This is the smallest possible value for~$\theta$.  If $\cos \theta = (a-c)/2$, then equation~(\ref{ident1}) gives $r_{12} = r_{14}$ and thus
$$
F \;  = \;  (r_{24}^3 - r_{14}^3)(r_{13}^3 - r_{12}^3)(r_{23}^3 - r_{34}^3) \;  >  \; 0,
$$ 
since $(a,b,c)$ is in the interior of~$D$.  (If any of the differences above also vanished, then $(a,b,c)$ would be on
the boundary of~$D$ due to Lemma~\ref{Lemma:Boundary}.)  Similarly, if 
$\cos \theta = (1-b)/(2a)$, then equation~(\ref{ident1}) gives $r_{12} = r_{23}$ and we compute that
\begin{eqnarray*}
F &  = &  (r_{13}^3 - r_{12}^3) \left[ (r_{24}^3 - r_{14}^3)(r_{12}^3 - r_{34}^3) -  (r_{12}^3 - r_{14}^3)(r_{24}^3 - r_{34}^3) \right]  \\[0.07in]
&  =  &     (r_{13}^3 - r_{12}^3) \left[ - r_{24}^3 r_{34}^3 - r_{14}^3 r_{12}^3  + r_{12}^3 r_{34}^3  + r_{14}^3 r_{24}^3 \right]  \\[0.07in]
&  =  &    (r_{13}^3 - r_{12}^3)(r_{24}^3 - r_{12}^3)(r_{14}^3 - r_{34}^3) \;  >  \; 0,
\end{eqnarray*}
since $(a,b,c)$ is in the interior of~$D$.
In either case, we see that $F(a,b,c, \theta = \theta_l) > 0$.

Next, suppose that  $\theta$ is chosen to be $\theta_u = \cos^{-1} (k_1)$.
This is the largest possible value for~$\theta$.  
If $\cos \theta = (c-a)/(2b)$, then equation~(\ref{ident2}) gives $r_{23} = r_{34}$ and thus
$$
F \;  = \;  -  (r_{12}^3 - r_{14}^3)(r_{24}^3 - r_{34}^3)(r_{13}^3 - r_{23}^3) \; < \;  0.
$$
If $\cos \theta = (b-1)/(2c)$, then equation~(\ref{ident2}) gives $r_{14} = r_{34}$ and we find that
\begin{eqnarray*}
F &  = &  (r_{24}^3 - r_{14}^3) \left[ (r_{13}^3 - r_{12}^3)(r_{23}^3 - r_{14}^3) -  (r_{12}^3 - r_{14}^3)(r_{13}^3 - r_{23}^3) \right]  \\[0.07in]
&  =  &     (r_{24}^3 - r_{14}^3) \left[ r_{13}^3 r_{23}^3 + r_{12}^3 r_{14}^3  - r_{12}^3 r_{13}^3  - r_{14}^3 r_{23}^3 \right]  \\[0.07in]
&  =  &    (r_{24}^3 - r_{14}^3)(r_{13}^3 - r_{14}^3)(r_{23}^3 - r_{12}^3) \;  <  \; 0,
\end{eqnarray*}
where the strict inequality follows once again from Lemma~\ref{Lemma:Boundary}.  If $\cos \theta = (a^2 - b^2 - 2b)/(2a)$, then 
equation~(\ref{ident3}) gives $r_{13} = r_{12}$ and thus
$$
F \;  = \;  -  (r_{12}^3 - r_{14}^3)(r_{24}^3 - r_{34}^3)(r_{13}^3 - r_{23}^3) \; < \;  0.
$$
Finally, if $\cos \theta = (1-c^2-2ac)/(2a)$, then equation~(\ref{ident3}) gives $r_{24} = r_{12}$ and we find that
\begin{eqnarray*}
F &  = &  (r_{12}^3 - r_{14}^3) \left[ (r_{13}^3 - r_{12}^3)(r_{23}^3 - r_{34}^3) -  (r_{12}^3 - r_{34}^3)(r_{13}^3 - r_{23}^3) \right]  \\[0.07in]
&  =  &     (r_{12}^3 - r_{14}^3) \left[ r_{13}^3 r_{23}^3 + r_{12}^3 r_{34}^3  - r_{12}^3 r_{13}^3  - r_{34}^3 r_{23}^3 \right]  \\[0.07in]
&  =  &    (r_{12}^3 - r_{14}^3)(r_{13}^3 - r_{34}^3)(r_{23}^3 - r_{12}^3) \;  <  \; 0.
\end{eqnarray*}
In all four cases, we find that $F(a,b,c, \theta = \theta_u) < 0$.  Since $F$ is a continuous function with opposite signs
at $\theta = \theta_l$ and $\theta = \theta_u$, the intermediate value theorem implies there exists an angle $\theta \in (\theta_l, \theta_u)$
such that $F(a,b,c,\theta) = 0$.

\vs

\nin {\bf (ii)  Uniqueness:}
To see that this solution is unique, we show that $\frac{ \partial F}{\partial \theta} < 0$ for any $(a,b,c)$ in the interior of~$D$ and
any $\theta \in (\theta_l, \theta_u)$.  From equations (\ref{eq:MutDist1}) and~(\ref{eq:MutDist2}), we have
$$
\frac{\partial r_{12}}{\partial \theta} = \frac{a \sin \theta}{r_{12}}, \quad
\frac{\partial r_{23}}{\partial \theta} = \frac{- ab \sin \theta}{r_{23}}, \quad
\frac{\partial r_{14}}{\partial \theta} = \frac{- c \sin \theta}{r_{14}}, \quad
\frac{\partial r_{34}}{\partial \theta} = \frac{bc \sin \theta}{r_{34}}, \quad \mbox{and} \quad
\frac{\partial r_{13}}{\partial \theta} = \frac{\partial r_{24}}{\partial \theta} = 0.
$$
Then we compute
$$
\frac{\partial F}{\partial \theta} \; = \;  -3 \sin \theta \left( a  r_{12} \, \alpha_1 + ab r_{23} \, \alpha_2 + c r_{14} \, \alpha_3 + bc r_{34} \, \alpha_4 \right),
$$
where
\begin{eqnarray}
\alpha_1  & = &  (r_{24}^3 - r_{14}^3)(r_{23}^3 - r_{34}^3) +  (r_{24}^3 - r_{34}^3)(r_{13}^3 - r_{23}^3), \nonumber  \\[0.07in]
\alpha_2  & = &  (r_{24}^3 - r_{14}^3)(r_{13}^3 - r_{12}^3) +  (r_{24}^3 - r_{34}^3)(r_{12}^3 - r_{14}^3),  \nonumber  \\[0.07in]
\alpha_3  & = &  (r_{24}^3 - r_{34}^3)(r_{13}^3 - r_{23}^3) -  (r_{13}^3 - r_{12}^3)(r_{23}^3 - r_{34}^3),  \label{eq:alpha3} \\[0.07in]
\alpha_4  & = &  (r_{24}^3 - r_{14}^3)(r_{13}^3 - r_{12}^3) -  (r_{12}^3 - r_{14}^3)(r_{13}^3 - r_{23}^3).  \nonumber
\end{eqnarray}
By~(\ref{geomCond}) and Lemma~\ref{Lemma:Boundary}, both $\alpha_1$ and $\alpha_2$ are strictly positive.
After adding and subtracting $r_{23}^6$ to~$\alpha_3$, we can rewrite that expression as
\begin{equation}
\alpha_3 \; = \;  (r_{24}^3 - r_{23}^3)(r_{13}^3 - r_{23}^3) +  (r_{12}^3 - r_{23}^3)(r_{23}^3 - r_{34}^3),
 \label{eq:alpha3Pos}
\end{equation}
which is also strictly positive on the interior of~$D$.  Finally, we find that
$$
\alpha_1 + \alpha_4 \; = \;  (r_{24}^3 - r_{14}^3)(r_{13}^3 - r_{12}^3 +  r_{23}^3 - r_{34}^3) +  (r_{13}^3 - r_{23}^3)(r_{24}^3 - r_{12}^3 + r_{14}^3 - r_{34}^3),
$$
which is strictly positive by~(\ref{geomCond}).  The conditions $a > c, 1 > b,$ and $r_{12} > r_{34}$, which are valid on the interior of~$D$,
combine to yield $a r_{12} > bc r_{34}$.   Then we have
$$
a r_{12} \, \alpha_1 + bc r_{34} \, \alpha_4 \; > \;  bc r_{34} \, \alpha_1 + bc r_{34} \, \alpha_4  \; = \;  bc r_{34} (\alpha_1 + \alpha_4) \; > \; 0 .
$$
This shows that $\frac{ \partial F}{\partial \theta} < 0$, which proves uniqueness.  

\vs

By the implicit function theorem, there exists a differentiable function $\theta = f(a,b,c)$ on the interior of~$D$ such that
$F(a,b,c,f(a,b,c)) = 0$.  The point $(a,b,c,\theta=f(a,b,c))$ describes a convex central configuration with positive masses.
Since $k_2 - k_1$ approaches zero as $(a,b,c)$ approaches the boundary of~$D$, we may extend
the function $f$ continuously to the boundary of~$D$, where it is defined as $\theta = \theta_l = \theta_u$.

Finally, if $(a,b,c) \not \in D$, then Lemma~\ref{Lemma:DistBds} shows that one of the mutual distance inequalities in~(\ref{geomCond}) will be violated.
For example, if $c > a$, then either $r_{12} < r_{14}$ or $r_{23} < r_{34}$.  Likewise, if $c \geq \frac{1}{a}(b^2 + 2b)$,
then either $r_{12} < r_{14}$ or $r_{13} \leq r_{12}$.   In any case, such a configuration, assuming it is central, will contain a negative or zero mass.
It follows that $D$ is precisely the domain of the implicitly defined function~$f$ and that the projection of $E$ into $abc$-space equals~$D$.
\enpf

%%%%%%%%%%%%%%%%%%
\subsection{Properties of the angle between the diagonals}
%%%%%%%%%%%%%%%%%%

Next we focus on the possible values of the angle $\theta$ between the two diagonals, showing that it is always between
$60^\circ$ and $90^\circ$.  Moreover, the value of $\theta$ increases as the radial variable $c$ increases.

\begin{lemma}
Suppose that $(a,b,c) \in D$ and $\theta = \pi/2$.  Then $r_{12}^3 + r_{34}^3 \geq  r_{14}^3 + r_{23}^3$.
\label{Lemma:Gpos}
\end{lemma}

\pf
When $\theta = \pi/2$, the formulas in (\ref{eq:MutDist1}) and~(\ref{eq:MutDist2}) reduce to
$
r_{12}^2 = a^2 + 1,  r_{14}^2 = c^2 + 1, r_{23}^2 = a^2 + b^2,$
and  $r_{34}^2 = b^2 + c^2$.  Define the function $G(a,b,c) = r_{12}^3 + r_{34}^3 -  r_{14}^3 - r_{23}^3$.
Note that $G(a,b,c=a) = 0$ since $r_{12} = r_{14}$ and $r_{23} = r_{34}$ when $c = a$ (a kite configuration).
We compute that
$$
\frac{ \partial G}{\partial c} \; = \;  3r_{34}^2 \frac{\partial r_{34}}{\partial c} - 3 r_{14}^2  \frac{\partial r_{14}}{\partial c}
\; = \;  3c (r_{34} - r_{14})  \; \leq \; 0,
$$
because $b \leq 1$ on~$D$.  Since $G(a,b,c=a) = 0$, it follows that $G(a,b,c < a) \geq 0$, as desired.
\enpf

\begin{theorem}
For a convex central configuration with positive masses, the angle $\theta$ between the two diagonals satisfies 
$\pi/3 < \theta \leq \pi/2$.  If $\theta = \pi/2$, the configuration must be a kite.
\label{Thm:angles}
\end{theorem}

\pf
We first show that $\theta \leq \pi/2$.  For any point $(a,b,c) \in D$, we have $r_{34} \leq r_{14}$ and $r_{23} \leq r_{12}$.
If $\theta = \pi/2$, we also have $r_{23}^3 - r_{34}^3 \leq r_{12}^3 - r_{14}^3$ by Lemma~\ref{Lemma:Gpos}.
Thus, when $\theta = \pi/2$, we have
\begin{equation}
(r_{24}^3 - r_{14}^3)(r_{13}^3 - r_{12}^3)(r_{23}^3 - r_{34}^3) \leq  (r_{24}^3 - r_{34}^3)(r_{13}^3 - r_{23}^3)(r_{12}^3 - r_{14}^3),
\label{ineq:thetaPi2}
\end{equation}
since all factors in~(\ref{ineq:thetaPi2}) are non-negative and each factor on the left-hand side
of the inequality is less than or equal to the corresponding factor on the right.
This shows that $F(a,b,c,\theta = \pi/2) \leq 0$.   From the proof of Theorem~\ref{Thm:Main}, 
$\partial F/\partial \theta < 0$ on the interior of $D \times [\theta_l, \theta_u]$.  Thus, for a fixed point $(a,b,c)$ in the interior of $D$, the unique solution
to $F(a,b,c,\theta) = 0$ must satisfy $\theta \leq \pi/2$.

Next, from~(\ref{ThetaBounds}), we have that $2 \cos \theta \leq a - c$ and $2a \cos \theta \leq 1 - b$.  We have just shown that $\cos \theta \geq 0$,
and since $b > 0$ and $c > 0$ on the interior of~$D$, we conclude that
\begin{equation}
2 \cos \theta \; <  \;  2 \cos \theta + c \;  \leq a \;  \leq \;  \frac{1 - b}{2 \cos \theta} \;  < \;  \frac{1}{2 \cos \theta} .
\label{ineq:theta}
\end{equation}
It follows that $\cos^2 \theta < 1/4$, which means $\theta > \pi/3$.

Finally, inequality~(\ref{ineq:thetaPi2}) is strict unless $r_{14} = r_{34}$ and $r_{12} = r_{23}$, or a factor on each side
of the inequality vanishes.   By Lemma~\ref{Lemma:Boundary}, this can only occur if $(a,b,c)$ lies on the
boundary of~$D$.  Thus, $F(a,b,c,\theta = \pi/2) < 0$ on the interior of~$D$.  
Since $\theta = \cos^{-1}(\frac{a-c}{2})$ or $\theta = \cos^{-1}(\frac{1-b}{2a})$ on the boundary of~$D$, we see that
a central configuration with $\theta = \pi/2$ must satisfy either $a = c$ or $b = 1$.  By (\ref{EqualCond1}) or~(\ref{EqualCond2}),
the configuration must be a kite.  
\enpf

\begin{remark}
\begin{enumerate}
\item  The fact that a convex central configuration with perpendicular diagonals must be a kite was proven earlier by the authors in~\cite{CCR}.

\item  If $\theta = \pi/3$, then all inequalities in (\ref{ineq:theta}) must become equalities.  This can only happen at the point $P_1 = (1,0,0)$,
a vertex of~$\overline{D}$ corresponding to an equilateral triangle configuration with bodies 3 and~4
coinciding ($r_{34} = 0$).  
\end{enumerate}
\end{remark}

Next we show that the value of $\theta$ increases as we move upwards (increasing in~$c$) through the domain~$D$.  We will need the following lemma.
Recall that $E$ is the set of four-body convex central configurations with positive masses in our particular coordinate system.

\begin{lemma}
Consider the following three quantities:
\begin{eqnarray*}
\beta_1 & = &  (r_{13}^3 - r_{12}^3)(r_{23}^3 - r_{34}^3)  - (r_{13}^3 - r_{23}^3)(r_{12}^3 - r_{14}^3),  \\[0.07in]
\beta_2 & = &  (r_{13}^3 - r_{23}^3)(r_{24}^3 - r_{34}^3)  - (r_{13}^3 - r_{12}^3)(r_{23}^3 - r_{34}^3),  \\[0.07in]
\beta_3 & = &  (r_{13}^3 - r_{23}^3)(r_{12}^3 - r_{14}^3)  - (r_{13}^3 - r_{12}^3)(r_{24}^3 - r_{14}^3).  
\end{eqnarray*}
Then, $\beta_1 \geq 0, \, \beta_2 > 0, \, \beta_3 < 0,$ and $\beta_2 + \beta_3 \geq 0$ for any configuration in~$E$.
\label{Lemma:beta-est}
\end{lemma}

\pf
Since we are working in~$E$, the equation $F = 0$ implies 
\begin{equation}
(r_{13}^3 - r_{23}^3)(r_{12}^3 - r_{14}^3) \; = \;   \frac{(r_{24}^3 - r_{14}^3)(r_{13}^3 - r_{12}^3)(r_{23}^3 - r_{34}^3)}{r_{24}^3 - r_{34}^3} \, .
\label{Eq:Fred}
\end{equation}
Then we have
$$
\beta_1 \; = \;   (r_{13}^3 - r_{12}^3)(r_{23}^3 - r_{34}^3)\left(1 - \frac{r_{24}^3 - r_{14}^3}{r_{24}^3 - r_{34}^3} \right) \; = \; 
\frac{(r_{13}^3 - r_{12}^3)(r_{23}^3 - r_{34}^3)(r_{14}^3 - r_{34}^3)}{r_{24}^3 - r_{34}^3} ,
$$
which is non-negative due to the inequalities in~(\ref{geomCond}).

Note that the quantity $\beta_2$ is identical to $\alpha_3$ used in the proof of Theorem~\ref{Thm:Main} (equation~(\ref{eq:alpha3})).
By equation~(\ref{eq:alpha3Pos}), we see that $\beta_2 > 0$ on~$E$.

Next, using equation~(\ref{Eq:Fred}), we have
$$
\beta_3 \; = \;   (r_{13}^3 - r_{12}^3)(r_{24}^3 - r_{14}^3)\left( \frac{r_{23}^3 - r_{34}^3}{r_{24}^3 - r_{34}^3} - 1 \right) \; = \; 
- \frac{(r_{13}^3 - r_{12}^3)(r_{24}^3 - r_{14}^3)(r_{24}^3 - r_{23}^3)}{r_{24}^3 - r_{34}^3} ,
$$
which is strictly negative due to the inequalities in~(\ref{geomCond}).

Finally, we compute that
\begin{eqnarray*}
\beta_2 + \beta_3  & = &   (r_{13}^3 - r_{23}^3)(r_{12}^3 + r_{24}^3 - r_{14}^3 - r_{34}^3)  -  (r_{13}^3 - r_{12}^3)(r_{23}^3 + r_{24}^3 - r_{14}^3 - r_{34}^3)  \\[0.07in]
& = &  r_{13}^3 ( r_{12}^3 - r_{23}^3)  - r_{23}^3 (r_{24}^3 - r_{14}^3 - r_{34}^3)  +  r_{12}^3 (r_{24}^3 - r_{14}^3 - r_{34}^3) \\[0.07in]
& = &   (r_{12}^3 - r_{23}^3)(r_{13}^3 + r_{24}^3 - r_{14}^3 - r_{34}^3),
\end{eqnarray*}
which is non-negative on~$E$.    This completes the proof.
\enpf

\begin{theorem}
On the interior of~$D$, the angle $\theta$ between the two diagonals increases with~$c$.  In other words, $\frac{\partial \theta}{\partial c} > 0$ on
the interior of~$D$.
\label{Thm:ThetaInc}
\end{theorem}

\pf
Recall that the angle $\theta = f(a,b,c)$ is a differentiable function on the interior of~$D$, determined by the solution to the equation
$F(a,b,c, f(a,b,c)) = 0$.  Using the implicit function theorem, we have
$
\frac{\partial \theta}{\partial c}  =    - \frac{\partial F}{\partial c} /  \frac{\partial F}{\partial \theta}   .
$
From the proof of Theorem~\ref{Thm:Main}, $\frac{\partial F}{\partial \theta} < 0$.   Thus, it suffices to show that $\frac{\partial F}{\partial c} > 0$, where
the partial derivative is evaluated at $(a,b,c, \theta = f(a,b,c)) \in E$ with $(a,b,c)$ in the interior of~$D$.

Using equations (\ref{eq:MutDist1}) and~(\ref{eq:MutDist2}), we have that
$$
\frac{\partial r_{14}}{\partial c} = \frac{c + \cos \theta}{r_{14}}, \quad
\frac{\partial r_{34}}{\partial c} = \frac{c  - b \cos \theta}{r_{34}}, \quad
\frac{\partial r_{24}}{\partial c} = 1, \quad  \mbox{and} \quad
\frac{\partial r_{12}}{\partial c} =  \frac{\partial r_{13}}{\partial c} = \frac{\partial r_{23}}{\partial c} = 0.
$$
Then we compute
\begin{equation}
\frac{\partial F}{\partial c} \;  =  \;  3r_{24}^2 \beta_1 + 3r_{14}(c + \cos \theta) \beta_2  + 3r_{34} (c - b \cos \theta) \beta_3 \, , 
\label{eq:PartialFc}
\end{equation}
where the $\beta_i$ are given as in Lemma~\ref{Lemma:beta-est}.   Since we are working in the interior of~$D$, the central configuration is not a kite and 
Theorem~\ref{Thm:angles} implies that $\cos \theta > 0$.  Hence, applying Lemma~\ref{Lemma:beta-est}, each term on the right-hand side of~(\ref{eq:PartialFc}) is non-negative except
for $3r_{34} c \beta_3$.  However, since $r_{14} \geq r_{34}$ and $\beta_2 > 0$, we have
\begin{eqnarray*}
3r_{14}c \beta_2 + 3r_{34} c \beta_3  & = &   3c(r_{14} \beta_2 + r_{34} \beta_3)  \\[0.05in]
& \geq &  3c( r_{34} \beta_2 + r_{34} \beta_3)  \\[0.05in]
& = & 3c r_{34} (\beta_2 + \beta_3) \\[0.05in]
& \geq & 0 
\end{eqnarray*}
by Lemma~\ref{Lemma:beta-est}.  This shows that $\frac{\partial F}{\partial c} > 0$.  The inequality is strict because the term $3r_{14} \beta_2 \cos \theta $
is strictly positive on the interior of $D$.  This completes the proof.  
\enpf

\begin{remark}
\begin{enumerate}
\item  Regarding Figure~\ref{Fig:setup}, if we fix the values of $a$ and $b$, then one interpretation of Theorem~\ref{Thm:ThetaInc} is that 
as the configuration widens in the vertical direction ($c$ increasing), the diagonals become closer and closer to perpendicular.   If $(a,b)$ is chosen from 
sub-region i or ii (see Figure~\ref{Fig:Dproj}), then the angle $\theta$ increases monotonically to $\pi/2$ where $c = a$ (a kite configuration).  On the other hand, if
$(a,b)$ belongs to sub-region iii or iv, then $\theta$ is bounded above by $\cos^{-1}(\frac{a-\bar{c}}{2}) < \pi/2$ where
$\overline{c} = \frac{1}{a}(b^2 + 2b) < a$.

\item  For kite configurations lying on the vertical face~II ($b = 1, r_{12} = r_{23},$ and $r_{14} = r_{34}$), it is straight forward to check that
$\frac{\partial F}{\partial c} = 0$.  This in turn implies that $\frac{\partial \theta}{\partial c} = 0$, which agrees with the fact that $\theta$
is constant ($\theta = \pi/2$) on all of face~II.  Thus, the strict inequality of Theorem~\ref{Thm:ThetaInc} only holds on the interior of~$D$.

\end{enumerate}
\end{remark}

%%%%%%%%%%%%%%%%%%%%%%%%%%
\section{Special Classes of Central Configurations}
\label{Sec:Shapes}
%%%%%%%%%%%%%%%%%%%%%%%%%%

In this section we use our coordinates in~$D$ to classify different types of quadrilaterals that are also central configurations.  
The analysis and defining equations are remarkably simple in our coordinate system, resulting in only linear or quadratic equations
in $a, b$, and~$c$.  Certain cases can be handled quickly due to the constraints on the mutual distances given
by~(\ref{geomCond}).  For example, the only parallelogram that can be a central configuration is a rhombus.
Likewise, the only possible rectangle is a square.

%%%%%%%%%%%%%%%%%%
\subsection{Kites}
%%%%%%%%%%%%%%%%%%

The kite configurations play a particularly important role in the overall classification of convex central configurations, 
occupying two of the six boundary faces of~$D$.  Recall that a kite configuration is a symmetric quadrilateral
with two bodies lying on the axis of symmetry and two bodies located equidistant from that axis.
The diagonals are always perpendicular and the two bodies not lying on the axis of symmetry must have equal mass.

\begin{figure}[htb]
 \centering
 \includegraphics[height=8.2cm,keepaspectratio=true]{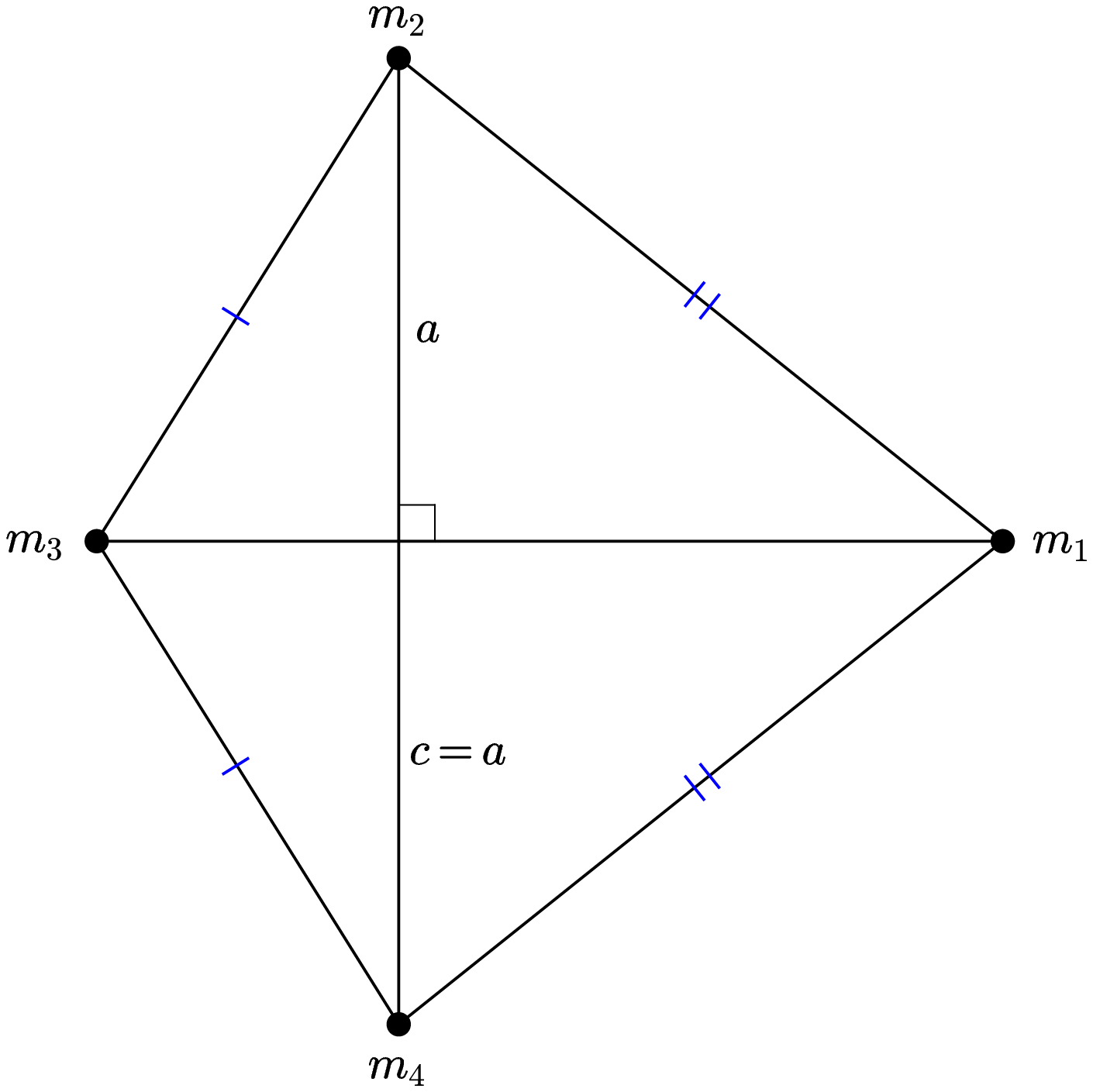}
 \hspace{-0.1in}
 \includegraphics[height=8.8cm,keepaspectratio=true]{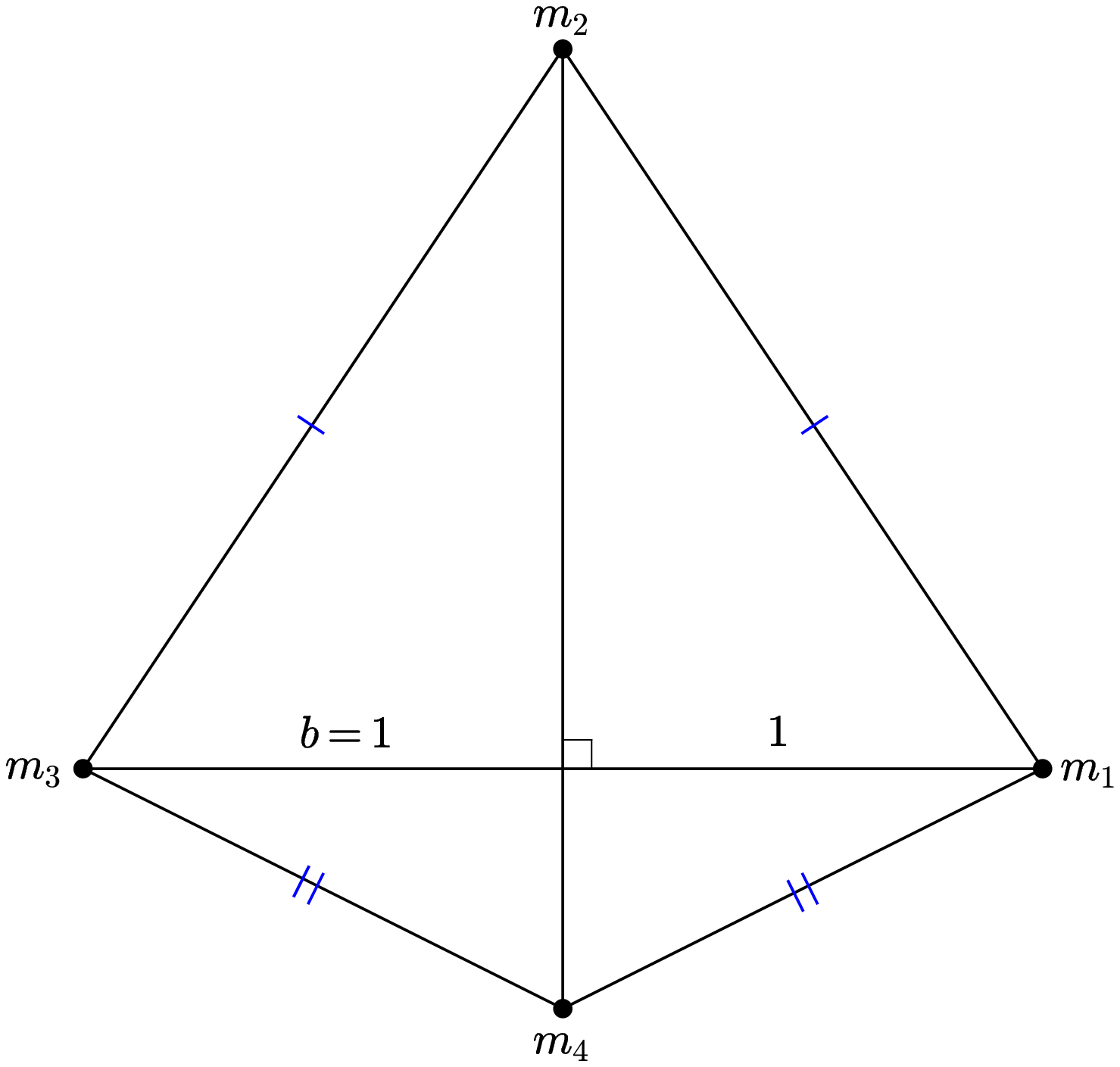}
 \caption{Two kite central configurations with different axes of symmetry.  Kites with a horizontal axis of symmetry (kite$_{13}$) 
 lie in the plane $c=a$, while those with a vertical axis of symmetry (kite$_{24}$) lie in the plane $b=1$.  All kites have $\theta = \pi/2$; these are the only possible
convex central configurations with perpendicular diagonals.}
 \label{Fig:kites} 
\end{figure}

Based on our ordering of the bodies, there are two possible types of kite configurations.
A kite with bodies 1 and~3 on the axis of symmetry, denoted {\em kite$_{13}$}, is symmetric with
respect to the $x$-axis and must satisfy $c = a$ (left plot in Figure~\ref{Fig:kites}).   These kites lie
on face~I and have $m_2 = m_4$, as can be verified by the middle formula in~(\ref{eq:massRatios2}).
A kite with bodies 2 and~4 on the axis of symmetry, denoted {\em kite$_{24}$}, is symmetric with
respect to the $y$-axis and must satisfy $b = 1$ (right plot in Figure~\ref{Fig:kites}).
These kites occupy face~II and require $m_1 = m_3$, 
as can be checked using the middle formula in~(\ref{eq:massRatios1}).

\vspace{-0.15in}

%----------------------------------------------------
%------------------------------------------------

It is important to note that due to statements (\ref{EqualCond1}) and~(\ref{EqualCond2}) in Lemma~\ref{Lemma:Boundary}, 
any point in $D$ lying on one of the two planes $c = a$ or $b = 1$ {\em must} correspond to a kite central configuration.
While two pairs of mutual distances must be congruent in order to distinguish a kite configuration from a general convex quadrilateral, only one equation
is required to imply a kite when restricting to the set of convex central configurations.   
An alternative interpretation of this fact is the following theorem.

\begin{theorem}
A convex central configuration with one diagonal bisecting the other must be a kite.
\label{thm:bisect}
\end{theorem}

\pf
In our coordinate system, if one of the diagonals bisects the other, then either $a = c$ or $b = 1$.  By (\ref{EqualCond1}) and~(\ref{EqualCond2}) in Lemma~\ref{Lemma:Boundary},
either case must correspond to a kite configuration.
\enpf

\begin{remark}
Theorem~\ref{thm:bisect} also follows directly from Conley's Perpendicular Bisector Theorem~\cite{rick2}.
\end{remark}

The intersection of the planes $c = a$ and $b = 1$ is a line that corresponds to the one-dimensional family
of rhombii central configurations.  This line is an edge on the boundary of~$D$ between vertices $P_3$ and $P_4$.
We regard $a$ as a parameter describing this family, with $1/\sqrt{3} < a < \sqrt{3}$.
From (\ref{eq:massRatios1}) and~(\ref{eq:massRatios2}), we have $m_1 = m_3, m_2 = m_4,$ and
$$
\frac{m_2}{m_1} \; = \;  \frac{8a^3 - a^3(a^2 + 1)^{3/2} }{ 8a^3 - (a^2 + 1)^{3/2}} \, .
$$
Note that $m_1$ and $m_3$ vanish as $a \rightarrow 1/\sqrt{3}$, while $m_2$ and $m_4$ approach 0 as 
$a \rightarrow \sqrt{3} \, $.   The length of the diagonal $r_{24}$ increases with $a$, stretching the rhombus in the vertical direction.  
The point $a = 1$ corresponds to the equal mass square configuration with congruent diagonals ($r_{13} = r_{24} = 2$).

%%%%%%%%%%%%%%%%%%
\subsection{Trapezoids}
%%%%%%%%%%%%%%%%%%

\begin{figure}[htb]
 \centering
 \includegraphics[height=8.8cm,keepaspectratio=true]{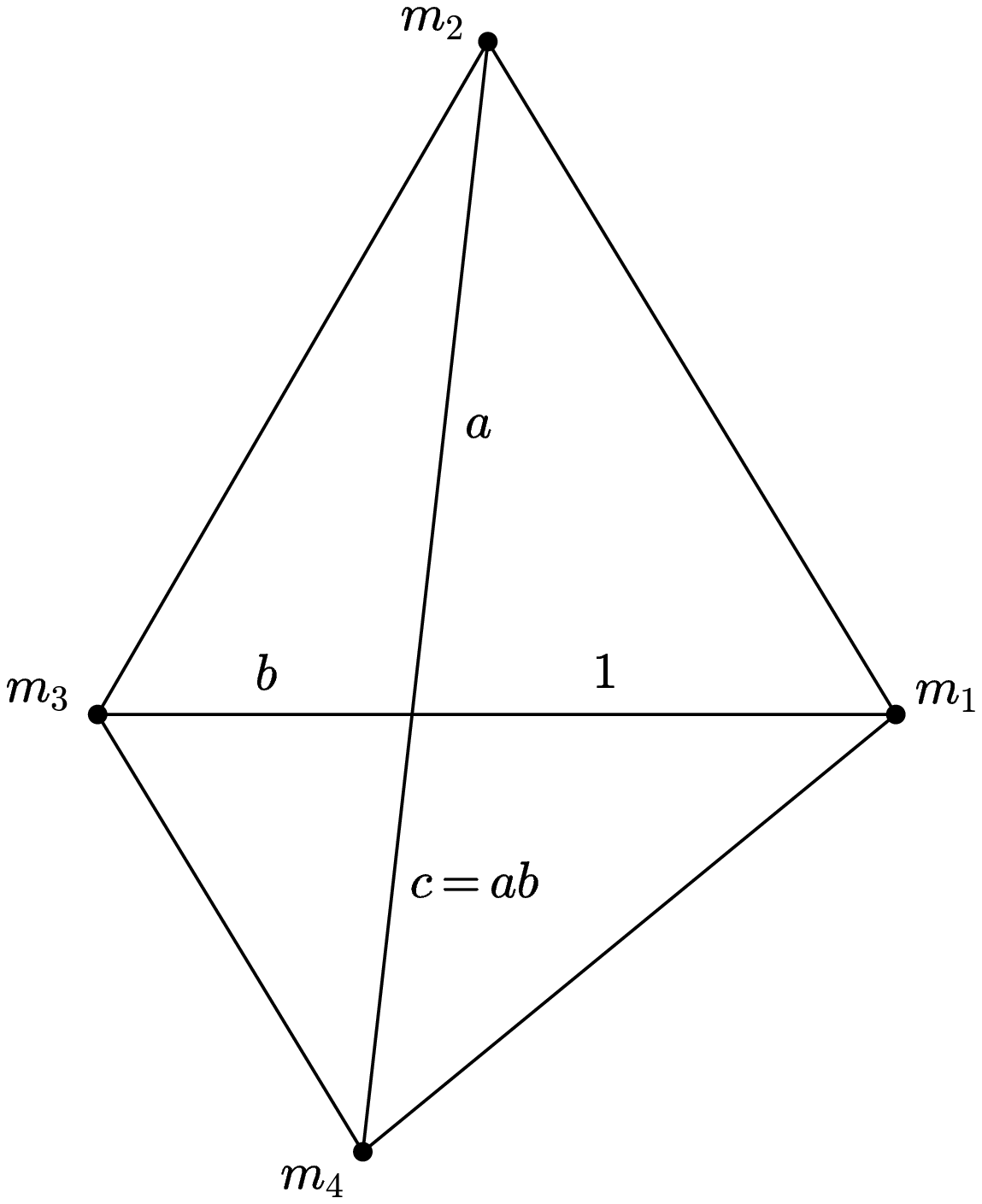}
 \hspace{0.05in}
 \includegraphics[height=8.2cm,keepaspectratio=true]{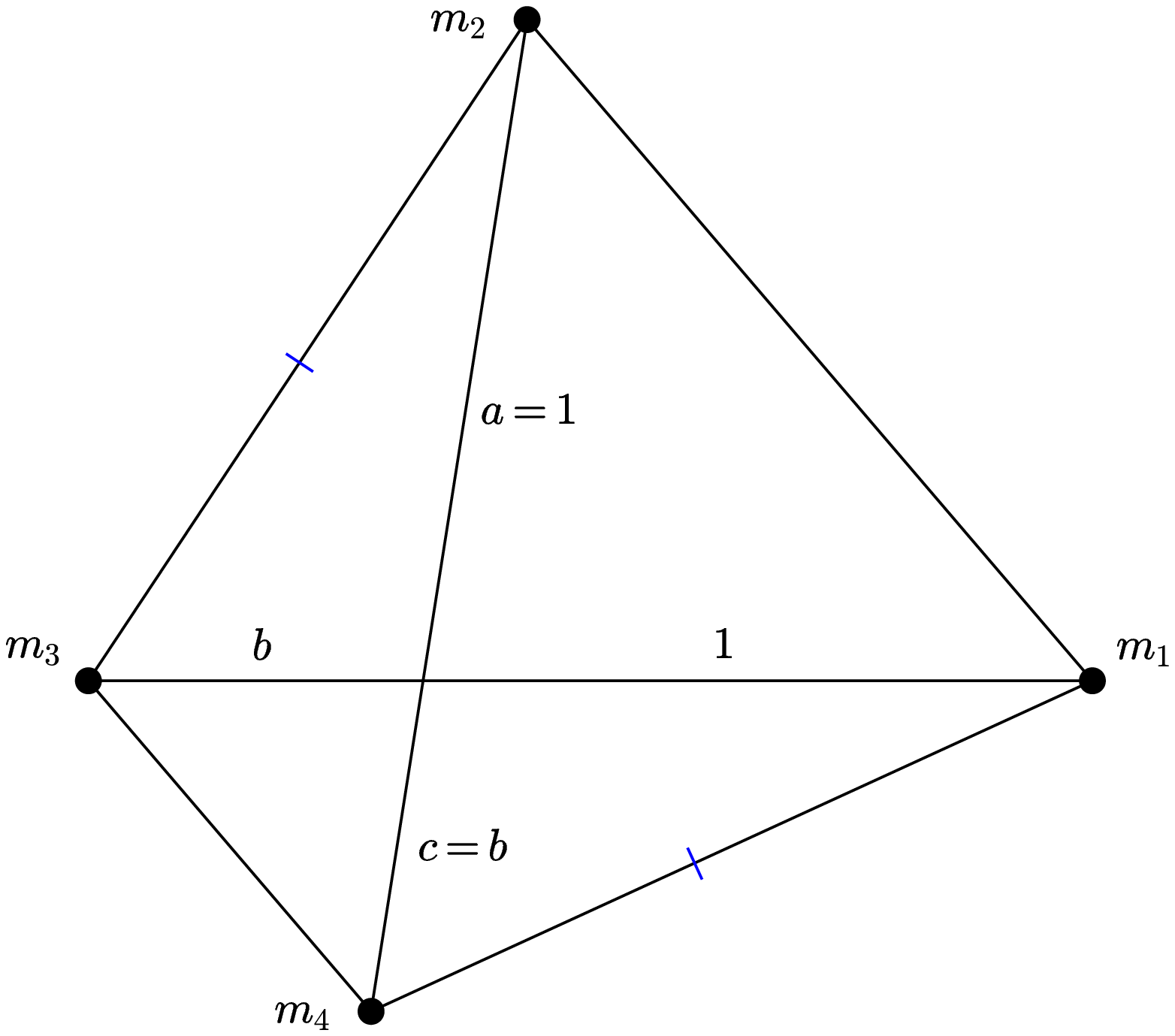}
 \vspace{-0.1in}
 \caption{Trapezoidal central configurations lie on the surface $c = ab$.  The isosceles trapezoid family (right figure)
 lies on the line formed by the intersection of the planes $a = 1$ and $c = b$.}
 \label{Fig:traps} 
\end{figure}

Next we consider the two possible types of trapezoids.  Let $\overline{q_i q_j}$ denote the side of the trapezoid
between vertices $i$ and~$j$.  If exterior sides $\overline{q_1 q_2}$ and $\overline{q_3 q_4}$ are parallel,
then we have
$$
\frac{a \sin \theta}{a \cos \theta - 1}  \; = \; \frac{c \sin \theta}{c \cos \theta - b} \, ,
$$  
which reduces to $(ab - c) \sin \theta =  0.$
Since $\sin \theta \neq 0$, $c = ab$ is both necessary and sufficient to have a trapezoid of this kind
(left plot in Figure~\ref{Fig:traps}).
On the other hand, if $\overline{q_1 q_4}$ is parallel to $\overline{q_2 q_3}$, 
then we quickly deduce that $a = bc$.  However, since
$a \geq c$ and $1 \geq b$ on~$D$, we have $a \geq bc$ always, with equality
only if both $a = c$ and $b =1$ are satisfied.  It follows that the only trapezoid
of this type is necessarily a rhombus, a subset of the first type of trapezoids.   This proves the following theorem.

\begin{theorem}
Suppose that $s$ is a central configuration in $E$.  Then $s = (a,b,c,\theta)$ is a trapezoid
if and only if $c = ab$.   The exterior sides $\overline{q_1 q_2}$ and $\overline{q_3 q_4}$ are always parallel.
\label{Thm:Trap}
\end{theorem}

%----------------------------------------------------
%-----------------------------------------------

\begin{remark}
Theorems \ref{Thm:Main} and~\ref{Thm:Trap} together show that the set of trapezoidal central configurations with positive masses is two-dimensional,
a graph over the surface $c = ab$ in~$D$ (a portion of a saddle).  This concurs with the recent results in~\cite{CC-Trap}.   \end{remark}

%----------------------------------------------------
 \begin{figure}[tb]
 \centering
 \includegraphics[height=9.5cm,keepaspectratio=true]{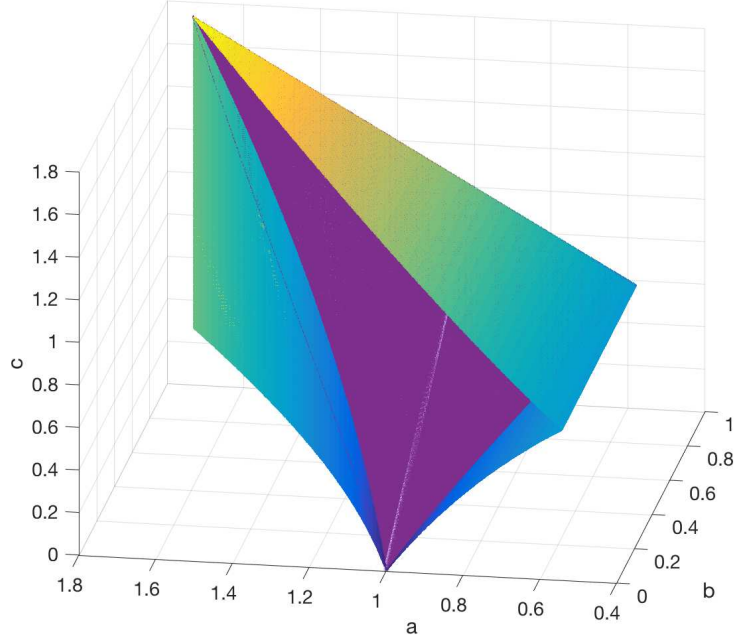}
 \vspace{-0.2in}
 \caption{The trapezoidal central configurations (purple) lie on the surface $c = ab$ within~$D$. 
 The violet line shows the isosceles trapezoid central configurations, where $a = 1$ and $c = b$.}
 \label{Fig:TrapsInD} 
\end{figure}
%-----------------------------------------------

Figure~\ref{Fig:TrapsInD} demonstrates how the surface of trapezoidal central configurations lies within the full space~$D$.
This surface intersects the boundary of~$D$ along the straight edge between vertices $P_3$ and~$P_4$ corresponding to the rhombii family (the intersection of faces I and~II).
It also meets the boundary of~$D$ in two curves of equilibria of the restricted four-body problem, one curve on face~V connecting vertices
$P_1$ and $P_4$, the other on face~IV joining vertices $P_1$ and~$P_3$.

Next, suppose that $s \in E$ is a trapezoid.  If we substitute $c = ab$ into equations (\ref{eq:MutDist1}) and~(\ref{eq:MutDist2}), we obtain
\begin{equation}
r_{23}^2 - r_{14}^2 \; = \;  (a^2 - 1)(1 - b^2)  \, .
\label{eq:trapIden}
\end{equation}
If $b = 1$, then $c = ab$ implies $c = a$ and hence $s$ is a rhombus.  Assuming that $b < 1$,  it follows from
equation~(\ref{eq:trapIden}) that $r_{23} > r_{14}$ for $a > 1$, and $r_{14} > r_{23}$ when $a < 1$.  The border between these
two cases are the isosceles trapezoids, where $r_{23} = r_{14}$ (right plot in Figure~\ref{Fig:traps}).  
In other words, the isosceles trapezoid family of central configurations corresponds to a line formed
by the intersection of the planes $a=1$ and $c = b$.  This line slices through the interior of $D$, crossing from
the degenerate equilateral triangle at $(1,0,0)$ to the square at $(1,1,1)$ (violet line in Figure~\ref{Fig:TrapsInD}).  By Theorem~\ref{Thm:ThetaInc}, the angle between
the diagonals monotonically increases from $\pi/3$ to~$\pi/2$ as $c$ increases from 0 to~1.  The family of isosceles trapezoids was studied
in \cite{pitu-gr} and \cite{xie}.

%%%%%%%%%%%%%%%%%%
\subsection{Co-circular configurations}
%%%%%%%%%%%%%%%%%%

\begin{figure}[t]
 \centering
 \includegraphics[height=7.5cm,keepaspectratio=true]{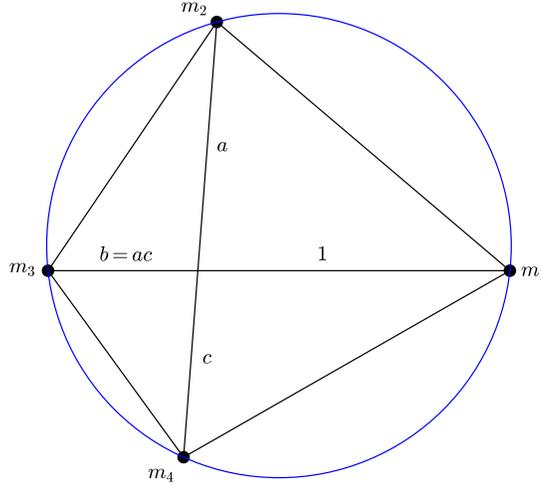}
 \vspace{-0.1in}
 \caption{A co-circular central configuration, where the bodies all lie on a common circle, must satisfy $b = ac$.}
 \label{Fig:co-circ} 
\end{figure}

Another interesting class of central configurations are those where the four bodies lie on a common circle, a {\em co-circular} central configuration 
(see Figure~\ref{Fig:co-circ}).
One of the main results in~\cite{pitu-gr} is that the set of four-body co-circular central configurations is a two-dimensional surface, a graph over
two of the exterior side-lengths.  We reproduce that result here, showing that the co-circular central configurations are a graph over the saddle $b = ac$ in~$D$.

%----------------------------------------------------
%-----------------------------------------------

\begin{theorem}
Suppose that $s$ is a central configuration in $E$.  Then $s = (a,b,c,\theta)$ is a co-circular central configuration 
if and only if $b = ac$. 
\end{theorem}

\pf
We make use of the cross ratio\footnote{Thanks to Richard Montgomery for suggesting this idea to the third author at the 2018 Joint Math Meetings.} 
from complex analysis~\cite{ahlfors}.  The cross ratio of four points $z_1, z_2, z_3, z_4$ is defined as the image of $z_1$ under the linear transformation that
maps $z_2$ to 1, $z_3$ to 0, and $z_4$ to $\infty$.  It is given by the expression
\begin{equation}
\frac{ (z_1 - z_3)(z_2 - z_4) }{ (z_1 - z_4)(z_2 - z_3) } .
\label{eq:cross-ratio}
\end{equation}
One of the nice properties of the cross ratio is that it is real if and only if the four points lie on a circle or a line.  Regarding the position of each body as a 
point in $\mathbb{C}$, we have $z_1 = 1, z_2 = ae^{i \theta}, z_3 = -b$, and $z_4 = -ce^{i \theta}$.  Substituting into~(\ref{eq:cross-ratio}), we find the
cross ratio to be
$$
\frac{ (a+c)(b+1) }{ ac e^{i \theta} + b e^{-i \theta} + a + bc} \, ,
$$
which is real if and only if $\sin \theta (ac - b) = 0$.  Since $\theta \in (\pi/3, \pi/2]$, we obtain $b = ac$ as a necessary and sufficient condition
for the four bodies to be lying on a common circle.
\enpf

\begin{figure}[t]
 \centering
 \includegraphics[height=9.5cm,keepaspectratio=true]{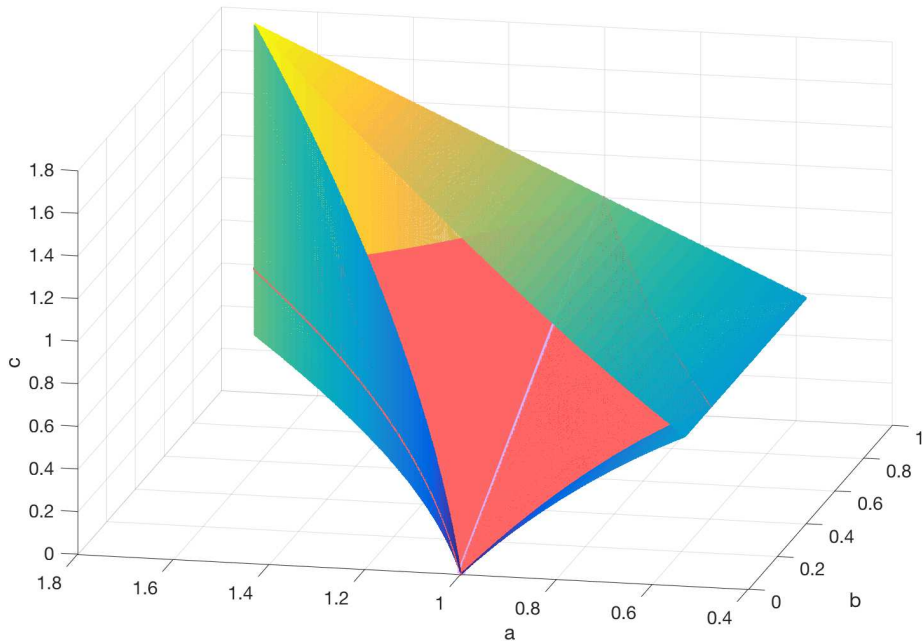}
 \vspace{-0.1in}
 \caption{Co-circular central configurations lie on the surface $b = ac$ (light red) in~$D$.  The violet line corresponds to the isosceles trapezoid family, where $a = 1$ and $b = c$.}
 \label{Fig:CoCircInD} 
\end{figure}

In Figure~\ref{Fig:CoCircInD} we plot the surface of co-circular central configurations within~$D$.  This surface intersects the boundary of~D on four faces.
On face~I we have co-circular kite configurations (kite$_{13}$) defined by the parabola $c = a, b = a^2, 1/\sqrt{3} < a \leq 1$.
We also have co-circular kites on face~II with the opposite axis of symmetry (kite$_{24}$), lying on the hyperbola $b = 1, c = 1/a, 1 \leq a < \sqrt{3}$.
The surface $b = ac$ also intersects faces IV and~V, tracing out curves of equilibria solutions to the restricted three-body problem.

Substituting $b = ac$ into equations (\ref{eq:MutDist1}) and~(\ref{eq:MutDist2}), we find that $r_{23} = a \, r_{14}$ and $r_{34} = c \, r_{12}$.  The line $a=1$ (violet line in
Figure~\ref{Fig:CoCircInD}) divides the surface $b = ac$ into two pieces.  As was the case for the trapezoids, if $1 < a < \sqrt{3}$, then we have co-circular central
configurations with $r_{23} > r_{14}$, while if $1/\sqrt{3} < a < 1$, then $r_{14} > r_{23}$.
Configurations on the line $a=1$ are isosceles trapezoids, where $r_{14} = r_{23}$.   Since $r_{12} \geq r_{34}$, the equation
$r_{34} = c \, r_{12}$ implies that $c \leq 1$ for any co-circular central configuration.  The maximum value of $c$ occurs at the square $a = b = c = 1$.

%----------------------------------------------------
%-----------------------------------------------

%%%%%%%%%%%%%%%%%%
\subsection{Equidiagonal configurations}
%%%%%%%%%%%%%%%%%%

\begin{figure}[h!]
 \centering
 \includegraphics[height=8.0cm,keepaspectratio=true]{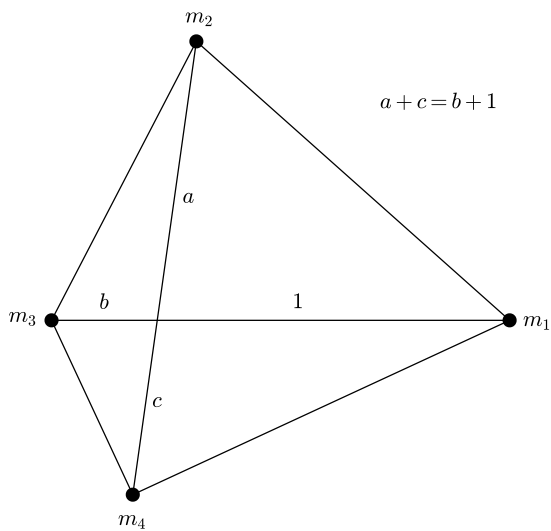}
 \hspace{-0.2in}
 \includegraphics[height=8.3cm,keepaspectratio=true]{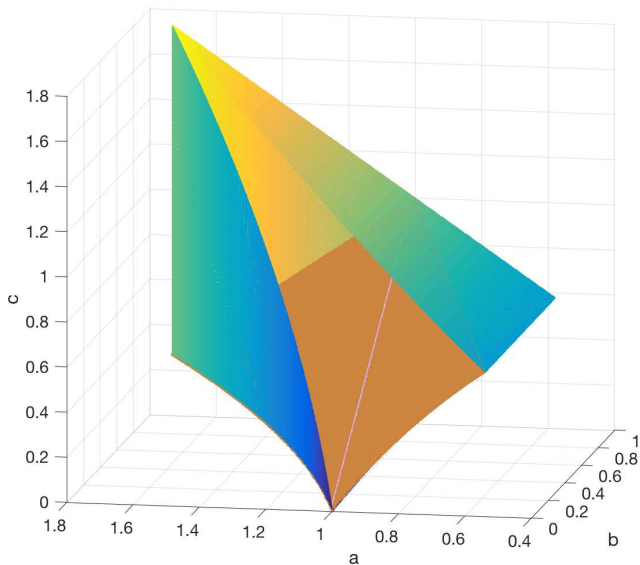}
 \vspace{-0.1in}
 \caption{Equidiagonal central configurations (left) are located on the plane $a - b + c = 1$ in~$D$ (right). 
 The violet line consists of the isosceles trapezoids ($a=1$ and $b=c$).}
 \label{Fig:EquiDiag} 
\end{figure}

The final class of convex central configurations we choose to explore are {\em equidiagonal} quadrilaterals, 
where the two diagonals are congruent (left plot in Figure~\ref{Fig:EquiDiag}).   These configurations are characterized by the equation
$r_{13} = r_{24}$, which is the plane $a - b + c = 1$ in our coordinates.  This plane intersects the boundary of~$D$ in four places
(right plot in Figure~\ref{Fig:EquiDiag}).  On face~I we find equidiagonal kites (kite$_{13}$) along the line
$c = a, b = 2a - 1, 1/\sqrt{3} < a \leq 1$.  Similarly, there is a line of equidiagonal kites (kite$_{24}$) on face~II parametrized by
$b = 1, c = 2 - a, 1 \leq a < \sqrt{3}$.  These two kite families  \linebreak

\vspace{-0.2in}

%----------------------------------------------------
%------------------------------------------------------

\noindent  intersect at the square $a = b = c = 1$.  
The equidiagonal plane also meets the boundary of~$D$ along two curved edges, 
one where faces III and~IV intersect and the other where faces V and~VI meet.  This follows directly from the equations given in
Table~\ref{Table:Faces}.

As with the trapezoidal and co-circular cases, the isosceles trapezoid family ($a = 1, b = c$) divides the equidiagonal plane into two regions
distinguished by whether $r_{23} > r_{14}$ (when $1 < a < \sqrt{3}$) or $r_{14} > r_{23}$ (when $1/\sqrt{3}< a < 1$).

%%%%%%%%%%%%%%%%%%
%\subsection{Tangential quadrilaterals}
%%%%%%%%%%%%%%%%%%

%A convex quadrilateral is {\em tangential} if the exterior sides are all tangent to a common circle.  Due to Pitot's theorem,
%a convex quadrilateral with our particular ordering of the bodies is tangential if and only if $r_{12} + r_{34} = r_{14} + r_{23}$.

%\begin{theorem}
%If a four-body convex central configuration is tangential, then it must be a kite.
%\label{Thm:tangential}
%\end{theorem}

%%%%%%%%%%%%%%%%%%
\subsection{Summary}
%%%%%%%%%%%%%%%%%%

Table~\ref{Table:Shapes} summarizes the different classes of configurations along with their defining equations
in $abc$-space or in the mutual distance variables~$r_{ij}$.  
In addition to the simplicity of the defining equations,
perhaps one of the more striking features of Table~\ref{Table:Shapes} is that
all of the configurations shown are defined by linear or quadratic equations.  
Moreover, due to Theorem~\ref{Thm:Main}, the dimension of each set is equivalent to the dimension of the 
corresponding geometric figure in $abc$-space.   Each type of configuration can be represented as the graph
of a function over a one- or two-dimensional set in~$D$, where the function is $\theta = f(a,b,c)$ restricted to
the given set.

\renewcommand{\arraystretch}{1.60}
\begin{table}[h!]
\begin{center}
\begin{tabular}{||c|c|c|c||}
\hline
\makebox[1.6in] {{\bf Configuration Type}} &  \makebox[1.4in]{{\bf Equation(s)}} & \makebox[1.8in]{{\bf Mutual Distances}} &\makebox[1.3in]{{\bf Figure in $D$}}  \\
\hline\hline
Kite$_{13}$ &  $c = a$   & $r_{12} = r_{14}$ and $r_{23} = r_{34}$  &   plane  \\
\hline
Kite$_{24}$ &  $b = 1$   & $r_{12} = r_{23}$ and $r_{14} = r_{34}$  &   plane  \\
\hline
Rhombus  &  $a = c$ and $b = 1$ & $r_{12} = r_{14} = r_{23} = r_{34}$  &  line \\
\hline
Trapezoid  &  $c = ab$ & $\overline{q_1 q_2}$ parallel to $\overline{q_3 q_4}$  & saddle  \\
\hline
Isosceles Trapezoid  &  $a = 1$ and $b = c$  & $r_{13} = r_{24}$ and $r_{14} = r_{23}$  &  line  \\
\hline
Co-circular  &   $b = ac$  &  $r_{13} r_{24} = r_{12} r_{34} + r_{14} r_{23}$  &  saddle  \\
\hline
Equidiagonal  &   $a - b + c = 1$   & $r_{13} = r_{24}$  &  plane \\
\hline\hline
\end{tabular}
\caption{Some special classes of central configurations and their surprisingly simple defining equations in $abc$-space.}
\label{Table:Shapes}
\end{center}
\end{table}
\renewcommand{\arraystretch}{1.0}

\begin{figure}[t]
 \centering
 \includegraphics[height=8.0cm,keepaspectratio=true]{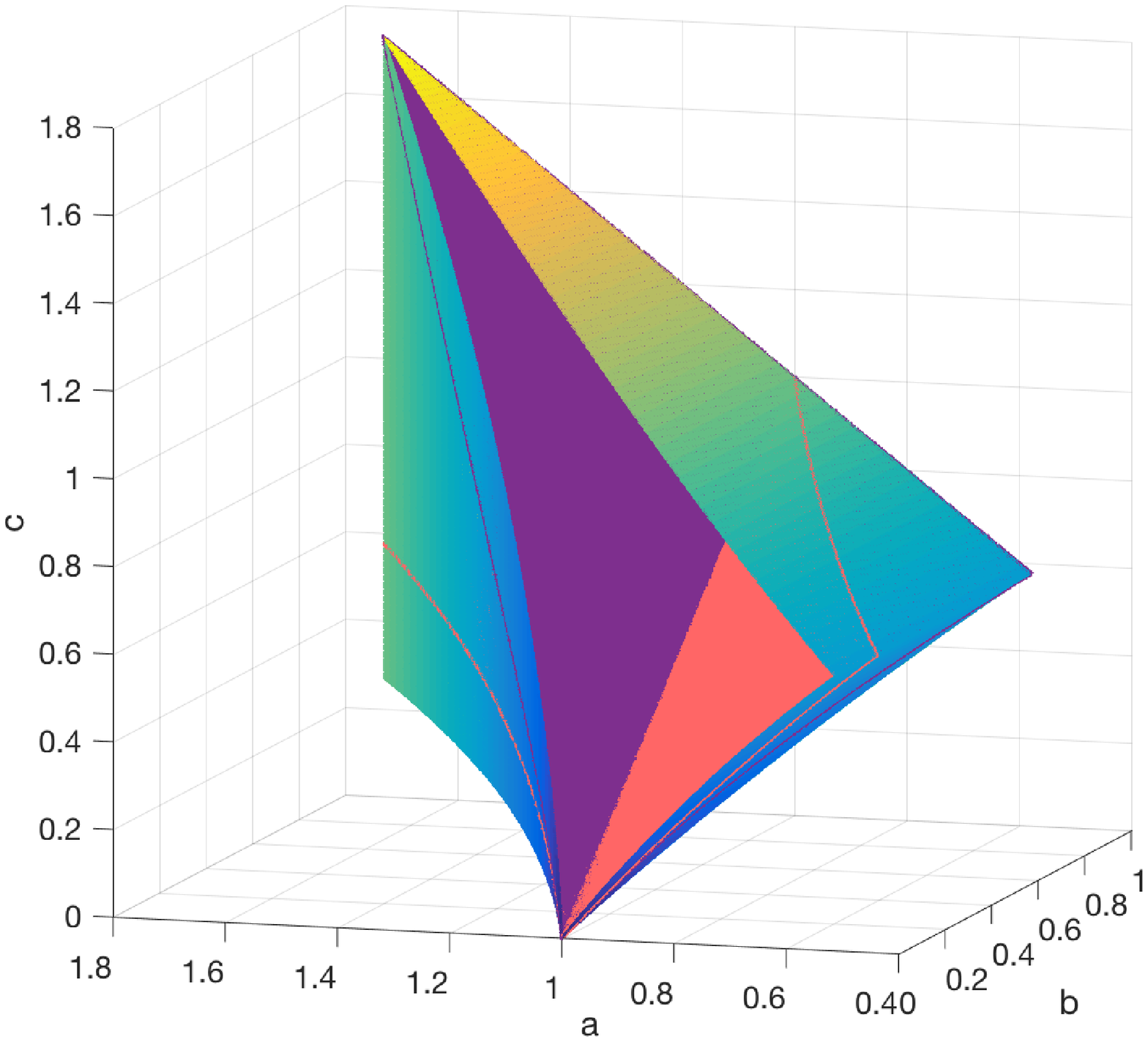}
 \hspace{-0.25in}
 \includegraphics[height=8.0cm,keepaspectratio=true]{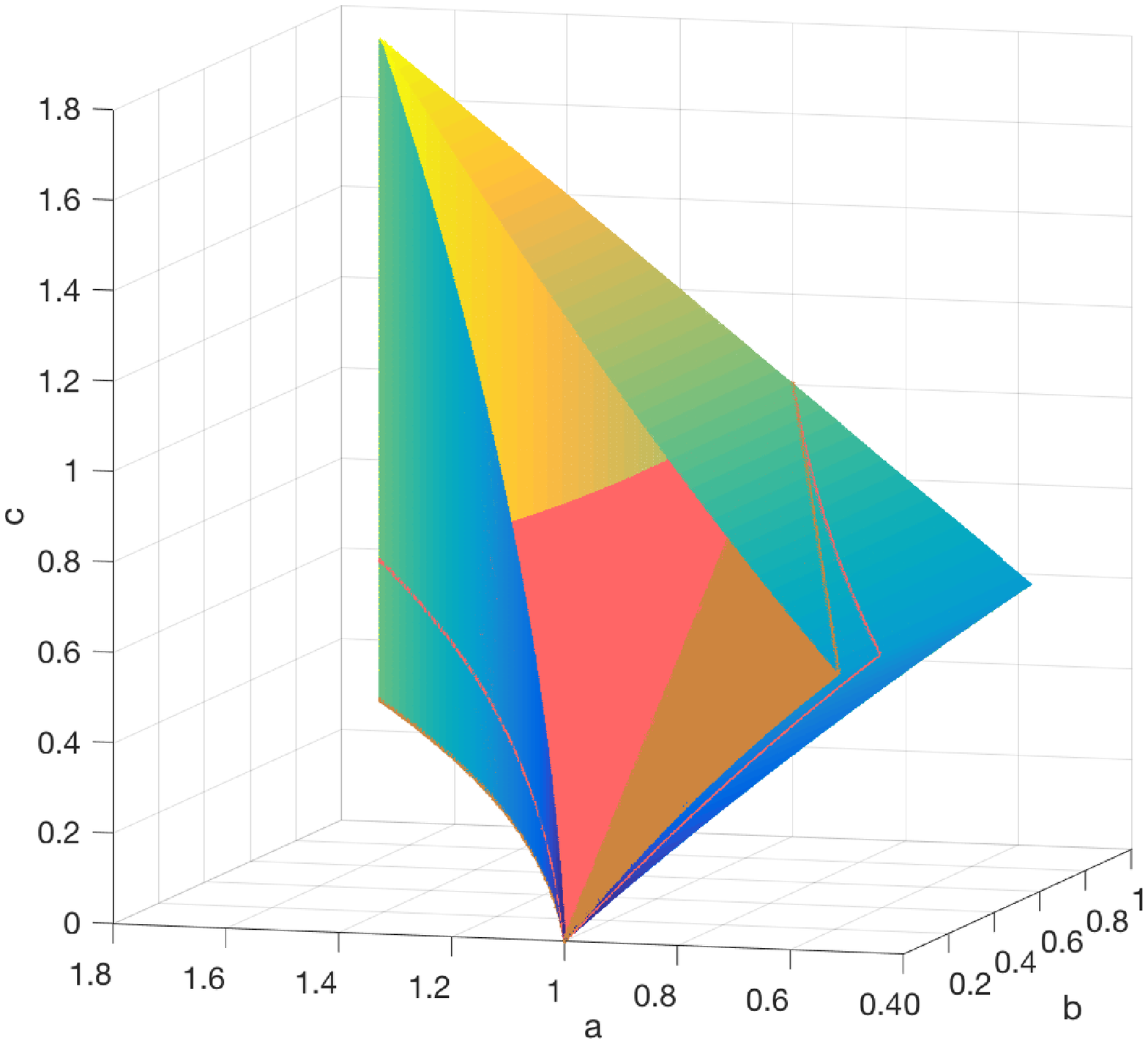}
 \caption{The left figure shows how the trapezoidal (purple) and co-circular (red) central configurations lie within~$D$, while 
 the right figure demonstrates how the co-circular (red) and equidiagonal (brown) central configurations fit together in~$D$.
 All three classes of central configurations intersect at the isosceles trapezoid family ($a = 1$ and $b = c$).}
 \label{Fig:MultSurfs} 
\end{figure}

Figure~\ref{Fig:MultSurfs} illustrates how the surfaces corresponding to trapezoidal, co-circular, and equidiagonal cnfigurations
lie within~$D$.   All three intersect at the line corresponding to the isosceles trapezoid configurations.  For $1 < a < \sqrt{3}$, the 
trapezoids are located above the co-circular configurations, which in turn lie above the equidiagonal solutions.  This is a consequence of
comparing the $c$-values on each surface.  Since $b \leq 1 < a$, we have
\begin{equation}
ab \; > \;  \frac{b}{a} \; > \; 1 - a + b \, .
\label{Ineq:3Surfs}
\end{equation}
On the other hand, for the portion of $D$ with $1/\sqrt{3} < a < 1$, the inequalities in~(\ref{Ineq:3Surfs}) are reversed and
the equidiagonal configurations lie above the co-circular solutions, which
lie above the trapezoids.

%----------------------------------------------------
%-------------------------------------------------------

The symmetric configurations play a particularly important role in the overall structure of~$D$, occupying two boundary faces
(kites), a boundary edge (rhombii), or a line of intersection between three classes of configurations (isosceles trapezoids).
Two classes of convex quadrilaterals must be kites in order to be central configurations.  
Configurations with either orthogonal or bisecting diagonals must be kites by Theorem~\ref{Thm:angles} and Theorem~\ref{thm:bisect}, respectively.

%%%%%%%%%%%%%%%%%%%%%%%
\section{Conclusion and Future Work}
%%%%%%%%%%%%%%%%%%%%%%%

We have established simple, yet effective coordinates for describing the space~$E$ of four-body convex central configurations.
Using these coordinates, we prove that $E$ is a three-dimensional set, the graph of a differentiable function 
over three radial variables.  The domain~$D$ of this function has been carefully defined, analyzed, and plotted
in $\mathbb{R}^{+^3}$.  Our coordinates provide elementary descriptions of several important classes of central configurations, including
kite, rhombus, trapezoidal, co-circular, and equidiagonal configurations.   The dimension and location of each of these classes within~$D$
has been explored in detail.  We have also shown that the angle between the diagonals of a four-body convex central configuration lies
between $60^\circ$ and $90^\circ$.  
As the configuration widens, the diagonals become closer and closer to orthogonal.  
The diagonals are perpendicular if and only if the quadrilateral is a kite.

In future research we intend to investigate the values of the masses as a function over the domain~$D$.
The mass ratios in (\ref{eq:massRatios1}) and~(\ref{eq:massRatios2}) reduce fairly nicely in our coordinate system,
although the dependence on the angle~$\theta = f(a,b,c)$ is complicated.
Nevertheless, we hope to build on our current work to show that the
mass map from $D$ into $\mathbb{R}^{+^3}$ (suitably normalized) is injective.  
Given a particular ordering of the bodies, this would prove that there
is a unique convex central configuration for any choice of four positive masses.

\vs

\noindent  {\bf Acknowledgments:}

M.~Corbera and J.~M.~Cors were partially supported by MINECO grant MTM2016-77278-P(FEDER);
J.~M.~Cors was also supported by AGAUR grant 2017 SGR 1617.
We also wish to thank John Little and Richard Montgomery for insightful discussions regarding this work.

%%%%%%%%%%%%%%%%%%%%%%%%%%%%%%%%%%%%%%%%%%%%%%%%%%%%%%%%%%%%%%%%%%%%%%%%%%%%%%%%%%
\bibliographystyle{amsplain}

\begin{thebibliography}{99}

\bibitem{albouy-symm1} Albouy, A., Sym\'{e}trie des configurations centrales de quatre corps, 
{\em C.~R. Acad. Sci. Paris}  {\bf 320}, s\'{e}rie 1 (1995), 217--220.


\bibitem{albouy-symm2}  Albouy, A., The symmetric central configurations of four equal masses, 
{\em Contemp. Math.}  {\bf 198} (1996), 131--135.


\bibitem{albouy-pblms}  Albouy, A., Cabral, H. E., and Santos, A. A.,
Some problems on the classical $n$-body problem,
{\em Celest. Mech. Dynam. Astronom.} {\bf 113} (2012), no. 4, 369--375.


\bibitem{albouy-chen}  Albouy, A. and Chenciner, A., 
Le probl\`{e}me des $n$~corps et les distances mutuelles,
{\em Invent. Math.}  {\bf 131} (1998), no. 1, 151--184. 


\bibitem{albouy} Albouy, A., Fu, Y., and Sun, S., 
Symmetry of planar four-body convex central configurations,
{\em Proc. R. Soc. Lond. Ser. A Math. Phys. Eng. Sci.}  {\bf 464} (2008), no. 2093, 1355--1365.


\bibitem{ahlfors}  Ahlfors, L. V., {\em Complex Analysis: An Introduction to the Theory of Analytic Functions of One Complex Variable}, 3rd ed., 
McGraw-Hill, New York, 1979.


%\bibitem{alb} Albouy, A., On a paper of Moeckel on central configurations,
%{\em Regul. Chaotic Dyn.} {\bf 8} (2003), no. 2, 133-142.

%\bibitem{barry} Barry, A. M., Hall, G. R. and Wayne, C. E.,
%Relative equilibria of the $(1+N)$-vortex problem, {\em J. Nonlinear Sci.}
%{\bf 22} (2012), 63-83.

%\bibitem{aref} Aref, H., Relative equilibria of point vortices and the fundamental theorem
%of algebra, {\em Proc. R. Soc. A} {\bf 467} (2011), 2168--2184.


%\bibitem{crystals} Aref, H., Newton, P. K., Stremler, M. A., Tokieda, T. and Vainchtein, D. L.,
%Vortex crystals, {\em Advances in Applied Mechanics} {\bf 39} (2003), 1--79.


\bibitem{leandro1} Barros, J. and Leandro, E.~S.~G.,
Bifurcations and enumeration of classes of relative equilibria in the planar restricted four-body problem,
{\em SIAM J. Math. Anal.}  {\bf 46} (2014), no. 2, 1185--1203. 


\bibitem{leandro2}  Barros, J. and Leandro, E.~S.~G.,
The set of degenerate central configurations in the planar restricted four-body problem, 
{\em SIAM J. Math. Anal.}  {\bf 43} (2011), no. 2, 634--661.


%\bibitem{barry} Barry, A., Hall, G. R. and Wayne, C. E., 
%Relative equilibria of the $(1+n)$-vortex problem,
%{\em J. Nonlinear Sci.} {\bf 22} (2012), 63--83.

%\bibitem{cas} Casasayas, J., Llibre, J and Nunes, A., Central configurations
%of the planar $1+N$ body problem, {\em Celest. Mech. Dynam. Astronom.}
%{\bf 60} (1994), 273--288.

%\bibitem{corb} Corbera, M., Cors, J. M. and Llibre, J.,
%On the central configurations of the planar $1+3$ body problem,
%{\em Celest. Mech. Dyn. Astr.}
%{\bf 109} (2011), 27-43.


\bibitem{CC-Bifur}  Corbera, M., Cors, J. M., Llibre, J., and Moeckel, R.,
Bifurcation of relative equilibria of the $(1+3)$-body problem,
{\em SIAM J. Math. Anal.}  {\bf 47} (2015), no. 2, 1377--1404. 


\bibitem{CC-Trap} Corbera, M., Cors, J. M., Llibre, J., and P\'{e}rez-Chavela, E., 
Trapezoid central configurations,  {\em Appl. Math. Comput.}  {\bf 346} (2019), 127--142. 


\bibitem{CCR} Corbera, M., Cors, J. M., and Roberts, G. E.,
A four-body convex central configuration with perpendicular diagonals is necessarily a kite,
{\em Qual. Theory Dyn. Syst.} {\bf 17} (2018), no. 2, 367--374.


%\bibitem{montse} Corbera, M., Delgado, J., and Llibre, J.,
%On the existence of central configurations of $p$ nested $n$-gons,
%{\em Qual. Theory Dyn. Syst.} {\bf 8} (2009), no. 2, 255--265.


%\bibitem{chenciner} Chenciner, A., Are there perverse choreographies? 
%{\em New advances in celestial mechanics and Hamiltonian systems},
%Kluwer/Plenum, New York (2004), 63--76.


\bibitem{pitu-gr}  Cors, J. M. and Roberts, G. E., 
Four-body co-circular central configurations, {\em Nonlinearity} {\bf 25} (2012), 343--370.


%\bibitem{davis} Davis, C., Wang, W., Chen, S. S.,  Chen, Y., Corbosiero, K., DeMaria, M., 
%Dudhia, J., Holland, G., Klemp, J., Michalakes, J., Reeves, H., Rotunno, R., Snyder, C.
%and Xiao, Q.,  
%Prediction of landfalling hurricanes with the Advanced Hurricane WRF Model,
%{\em Monthly Weather Review} {\bf 136} (2007), 1990--2005.


%\bibitem{little} Cox, D., Little, J. and O'Shea, D., \textit{Using Algebraic Geometry},
%2nd. ed., Springer, New York, 2005.


\bibitem{dzio} Dziobek, O., \"{U}ber einen merkw\"{u}rdigen Fall des
Vielk\"{o}rperproblems, {\em Astro. Nach.} {\bf 152} (1900), 32--46.


\bibitem{erdi}  \'{E}rdi, B. and Czirj\'{a}k, Z.,
Central configurations of four bodies with an axis of symmetry,
{\em Celestial Mech. Dynam. Astronom.}  {\bf 125} (2016), no. 1, 33--70. 




\bibitem{fernand} Fernandes, A. C., Llibre, J., and Mello, L. F., 
Convex central configurations of the four-body problem with 
two pairs of equal adjacent masses, {\em Arch. Ration. Mech. Anal.} {\bf 226} (2017), no.~1, 303--320.


\bibitem{hall} Hall, G. R., Central configurations in the planar $1+n$ body problem, preprint (1988).



\bibitem{hampton}  Hampton, M., Concave central configurations in the four-body problem,
{\em Doctoral Thesis}, University of Washington, Seattle (2002).


\bibitem{rick-hamp}  Hampton, M. and Moeckel, R.,  Finiteness of relative equilibria of the four-body problem, 
{\em Invent. Math.} {\bf 163} (2006), 289--312.



%\bibitem{hamp}  Hampton, M., Co-circular central configurations in the four-body
%problem, {\em EQUADIFF 2003} (Conference Proceedings), World Sci. Publ., 
%Hackensack, NJ, (2005), 993--998.


\bibitem{HRS} Hampton, M., Roberts, G. E., and Santoprete, M., 
Relative equilibria in the four-vortex problem with two pairs of
equal vorticities,  {\em J. Nonlinear Sci.} {\bf 24} (2014), 39--92.


\bibitem{kule} Kulevich, J. L.,  Roberts, G. E., and Smith, C. J.,
Finiteness in the planar restricted four-body problem,
{\em Qual. Theory Dyn. Syst.} {\bf 8} (2009), no. 2, 357--370. 


\bibitem{lagrange}  Lagrange, J. L., Essai sur le probl\'{e}me des trois corps.,  {\em \OE uvres} {\bf 6} (1772),
Gauthier-Villars, Paris, 272--292.


%\bibitem{Li} Li, B., Deng, Y. and Zhang, S., A property for four-body isosceles trapezoid central configurations,
%preprint, \url{http://www.paper.edu.cn/download/downPaper/201606-641}


%\bibitem{Llibre} Llibre, J. and Valls, C., The co-circular central configurations of the 5-body problem, preprint, (April, 2013).


\bibitem{Mac} MacMillan, W. D. and Bartky, W., Permanent configurations in the problem of four bodies,
{\em Trans. Amer. Math. Soc.} {\bf 34} (1932), no. 4, 838--875.



%\bibitem{mars} Marsden, J. E. and Ross, S. D., New methods in celestial mechanics and mission design,
%{\it Bull. Amer. Math. Soc. (N.S.)} {\bf 43} (2006), no. 1, 43--73.

\bibitem{matlab}  MATLAB, version R2016b (9.1.0.441655), (2016), The MathWorks, Inc., 
Natick, Massachusetts, United States.


\bibitem{meyer} Meyer, K. R. and Offin, D. C., {\em Introduction to Hamiltonian Dynamical Systems
and the $N$-Body Problem}, 3rd ed., Applied Mathematical Sciences, 90, Springer, Cham, 2017.


%\bibitem{maple}  Maple, version 15.00, Maplesoft, Waterloo Maple Inc. (2011).


\bibitem{rick-book}  Moeckel, R., Central configurations, in
{\em Central Configurations, Periodic Orbits, and Hamiltonian Systems},
Llibre, J., Moeckel, R. and Sim\'{o}, C, Birkh\"{a}user (2015), 105--167.


\bibitem{rick1} Moeckel, R., Linear stability of relative equilibria with a dominant mass,
{\em Differ. Equ. Dyn. Syst.} {\bf 6} (1994), no. 1, 37--51.


\bibitem{rick2} Moeckel, R., On central configurations, {\em Mathematische Zeitschrift}
{\bf 205} (1990), no. 4, 499--517. 


\bibitem{rick-cluster}  Moeckel, R., Relative equilibria with clusters of small masses,
{\em J. Dynam. Differ. Equ.} {\bf 9} (1997), no. 4, 507--533.


%\bibitem{moulton}  Moulton, F. R.,  The straight line solutions of the problem of $n$~bodies, 
%{\em Ann. of Math. (2)} {\bf 12}  (1910), no. 1, 1--17. 


%\bibitem{newton} Newton, P. K., {\em The $N$-Vortex Problem: Analytic Techniques}, 
%Springer, New York, 2001.


%\bibitem{oneil}  O'Neil, K. A.,  Stationary configurations of point vortices, 
%{\em Trans. Amer. Math. Soc.}  {\bf 302}, no. 2 (1987),  383--425. 


%\bibitem{g:thesis}  Roberts, G. E., Existence and stability of relative equilibria in the $n$-body problem,
%{\em Doctoral Thesis}, Boston University (1999).


\bibitem{saari} Saari, D. G., {\em Collisions, Rings, and Other Newtonian $N$-Body Problems},
CBMS Regional Conference Series in Mathematics, no. 104, Amer. Math. Soc., Providence, RI, 2005.


\bibitem{sage} SageMath, the Sage Mathematics Software System (Version 7.3), The Sage Developers,
2016, \url{http://www.sagemath.org}.


\bibitem{santoprete}  Santoprete, M.,  Four-body central configurations with one pair of opposite sides parallel, 
{\em J. Math. Anal. Appl.}  {\bf 464} (2018) 421--434.



\bibitem{schmidt} Schmidt, D., Central configurations and relative
equilibria for the $n$-body problem,  {\em Classical and celestial
mechanics (Recife, 1993/1999)}, Princeton Univ. Press, Princeton, NJ (2002), 1--33.


\bibitem{simo}  Sim\'{o}, C., Relative equilibrium solutions in the four-body problem, 
{\em Celestial Mech.}  {\bf 18} (1978), no. 2, 165--184.
 

%\bibitem{smale} Smale, S., Mathematical problems for the next century,  
%{Math. Intelligencer} {\bf 20} (1998), no. 2, 7-15.


%\bibitem{sturm} Sturmfels, B., {\em Solving Systems of Polynomial Equations},
%CBMS Regional Conference Series in Mathematics, no. 97, 
%Amer. Math. Soc., Providence, RI, 2002.


\bibitem{wintner} Wintner, A., \textit{The Analytical Foundations of Celestial Mechanics},
Princeton Math. Series {\bf 5}, Princeton University Press, Princeton, NJ, 1941.



\bibitem{Xia} Xia, Z., Convex central configurations for the $n$-body problem,
{\em J. Diff. Equa.} {\bf 200} (2004), 185--190.



\bibitem{xie}  Xie, Z., Isosceles trapezoid central configurations of the Newtonian four-body problem,
{\em Proc. Roy. Soc. Edinburgh Sect. A} {\bf 142} (2012), no. 3, 665--672.




\end{thebibliography}

\end{document}